\documentclass[12pt,reqno]{amsart}

\title[Convex minimization problems]
{Convex minimization problems\\ with weak constraint
qualifications}
\author{Christian L\'eonard}
\date{September 07}

\usepackage{amssymb, amsmath, amsfonts, latexsym, enumerate}

\newtheorem{theorem}[equation]{Theorem}
\newtheorem{lemma}[equation]{Lemma}
\newtheorem{proposition}[equation]{Proposition}
\newtheorem{corollary}[equation]{Corollary}

\newtheorem{definition}[equation]{Definition}

\theoremstyle{remark}
\newtheorem{remark}[equation]{Remark}

\newtheorem{questions}[equation]{Questions}

\numberwithin{equation}{section}

\newcommand{\R}{\mathbb{R}}

\newcommand{\eqdef}{\stackrel{\vartriangle}{=}}
\newcommand{\dom}{\mathrm{dom\,}}
\newcommand{\icordom}{\mathrm{icordom\,}}

\newcommand{\inter}{\mathrm{int\,}}
\newcommand{\epi}{\mathrm{epi}\,}
\newcommand{\lsc}{lower semicontinuous}
\newcommand{\usc}{upper semicontinuous}
\newcommand{\cav}{{\hat{*}}}

\newcommand{\ls}{\mathrm{ls}\,}

\newcommand{\boulette}[1]{$\bullet$\ Proof of (#1).}
\newcommand{\Boulette}[1]{\par\medskip\noindent $\bullet$\ Proof of (#1).}

\newcommand{\diffdom}{\mathrm{diffdom\,}}

\newcommand{\UU}{\mathcal{U}}
\newcommand{\LL}{\mathcal{L}}
\newcommand{\YY}{\mathcal{Y}}
\newcommand{\XX}{\mathcal{X}}
\newcommand{\ZZ}{\mathcal{Z}}
\newcommand{\UUo}{\mathcal{U}_o}
\newcommand{\LLo}{\mathcal{L}_o}
\newcommand{\YYo}{\mathcal{Y}_o}
\newcommand{\XXo}{\mathcal{X}_o}
\newcommand{\Co}{C_o}
\newcommand{\Ci}{C}
\newcommand{\Ss}{S^{\ast}}
\newcommand{\Ct}{C_\theta}
\newcommand{\Cts}{C_{\theta^*}}

\newcommand{\AB}{{A\!\times\! B}}
\newcommand{\PA}{P_A}
\newcommand{\PB}{P_B}
\newcommand{\PAB}{P_{AB}}
\newcommand{\CAB}{C_{AB}}

\newcommand{\Fo}{\Phi_o}
\newcommand{\Fos}{\Phi_o^*}
\newcommand{\F}{\Phi}
\newcommand{\Fs}{\Phi^*}
\newcommand{\Fb}{\overline{\Phi}}
\newcommand{\Lao}{\Lambda_o}
\newcommand{\Los}{\Lambda_o^*}
\newcommand{\La}{\Lambda}
\newcommand{\Ls}{\Lambda^*}
\newcommand{\Lb}{\overline{\Lambda}}
\newcommand{\To}{T_o}
\newcommand{\fss}{\varphi^{**}}
\newcommand{\ts}{\theta^*}
\newcommand{\jt}{j_{\theta}}
\newcommand{\jts}{j_{\theta^*}}
\newcommand{\ds}{\iota_{\Ct}^*}

\newcommand{\Po}{$(P_o)$}
\renewcommand{\Pi}{$(P)$}
\newcommand{\PX}{$(P^{\XX})$}
\newcommand{\Pxh}{$(P^{{\hat x}})$}
\newcommand{\Do}{$(D_o)$}
\newcommand{\Di}{$(D)$}
\newcommand{\Dii}{$(\overline{D})$}
\newcommand{\Dxh}{$(\overline{D}^{{\hat x}})$}

\newcommand{\NF}{|\cdot|_\F}
\newcommand{\NL}{|\cdot|_\La}
\newcommand{\sLU}{\sigma(\LLo,\UUo)}
\newcommand{\sLUi}{\sigma(\LL,\UU)}
\newcommand{\sUL}{\sigma(\UUo,\LLo)}
\newcommand{\sULi}{\sigma(\UU,\LL)}
\newcommand{\sULii}{\sigma(\LL^\ast,\LL)}
\newcommand{\sXY}{\sigma(\XXo,\YYo)}
\newcommand{\sXYi}{\sigma(\XX,\YY)}
\newcommand{\sXYii}{\sigma(\XX,\XX^\ast)}
\newcommand{\sYXii}{\sigma(\XX^\ast,\XX)}
\newcommand{\sYX}{\sigma(\YYo,\XXo)}
\newcommand{\sYXi}{\sigma(\YY,\XX)}

\newcommand{\HF}{$(H_\F)$}
\newcommand{\HFi}{$(H_{\F1})$}
\newcommand{\HFii}{$(H_{\F2})$}
\newcommand{\HFiii}{$(H_{\F3})$}
\newcommand{\HT}{$(H_T)$}
\newcommand{\HTi}{$(H_{T1})$}
\newcommand{\HTii}{$(H_{T2})$}
\newcommand{\HC}{$(H_{C})$}

\newcommand{\ul}{\langle u,\ell\rangle}
\newcommand{\lz}{\langle \ell,\zeta\rangle}
\newcommand{\xo}{\langle x,\omega\rangle}
\newcommand{\yx}{\langle y,x\rangle}

\newcommand{\lh}{\hat\ell}
\newcommand{\ob}{\bar\omega}
\newcommand{\yb}{\bar y}
\newcommand{\xh}{\hat{x}}
\newcommand{\bq}{\langle b,q\rangle}
\newcommand{\ap}{\langle a,p\rangle}
\newcommand{\ab}{\bar a}
\newcommand{\bb}{\bar b}

\newcommand{\IAB}{\int_{\AB}}
\newcommand{\IZ}{\int_{\ZZ}}

\newcommand{\lmax}{{\lambda_\diamond}}

\newcommand{\MZ}{M_\ZZ}

\newcommand{\Ml}{M^{\lmax}_\ZZ}

\newcommand{\EEl}{\mathcal{E}_{\lmax}}

\newcommand{\remMK}{4.9} 

\begin{document}


 \address{Modal-X, Universit\'e Paris 10\quad \& }
 \address{CMAP, \'Ecole Polytechnique. 91128 Palaiseau Cedex, France}
 \email{christian.leonard@polytechnique.fr}

\begin{abstract}
One revisits the standard saddle-point method based on conjugate
duality for solving convex minimization problems. Our aim is to
reduce or remove unnecessary topological restrictions on the
constraint set. Dual equalities and characterizations of the
minimizers are obtained with weak or without constraint
qualifications.
\\
The main idea is to work with intrinsic topologies which reflect
some geometry of the objective function.
\\
The abstract results of this article are applied in other papers
to the Monge-Kantorovich optimal transport problem and the
minimization of entropy functionals.
\end{abstract}

\maketitle
\tableofcontents


\section{Introduction}\label{sec:introduction}

An ``\emph{extension}'' of the saddle-point method for solving a
convex minimization problem is investigated. It is shown how to
implement the standard saddle-point method in such a way that
topological restrictions on the constraint sets (the constraint
qualifications) may essentially be removed. With this aim in view,
one works with topologies associated with gauge functionals of
sets which are close to the level sets of the objective function.
These well-suited topologies partly reflect the geometry of the
problem.
\\
At some point, one has to compute an \emph{extended} dual problem.
This is the price to pay for implementing this approach.
\\
The method is based on conjugate duality as developed by
R.T.~Rockafellar in \cite{Roc74}. Dual equalities and
characterizations of the minimizers are obtained with weak or
without constraint qualification.
\\
This paper is a companion of \cite{p-Leo07b} and \cite{p-Leo07c}
where this extended saddle-point method is applied to the
Monge-Kantorovich optimal transport problem and the minimization
of entropy functionals.

\subsection*{An abstract convex problem and related questions}

Let $\UU$ be a vector space, $\LL=\UU^\ast$ its algebraic dual
space, $\F$ a $(-\infty,+\infty]$-valued convex function on $\UU$
and $\Fs$ its convex conjugate for the duality $\langle
\UU,\LL\rangle.$ Let $\YY$ be another vector space, $\XX=\YY^\ast$
its algebraic dual space and $T:\LL\rightarrow\XX$ is a linear
operator. We consider the convex minimization problem
\begin{equation}
  \textsl{minimize } \Fs(\ell) \textsl{ subject to } T\ell\in C,\
  \ell\in\LL
\tag{$P$}
\end{equation}
where $ C$ is a convex subset of $\XX.$ As is well known,
Fenchel's duality leads to the dual problem
\begin{equation}
  \textsl{maximize } \inf_{x\in C}\yx - \F(T^\ast y),\quad y\in\YY
\tag{$D$}
\end{equation}
where $T^\ast $ is the adjoint of $T.$
\begin{questions}\label{questions}
The usual questions related to $(P)$ and $(D)$ are
\begin{itemize}
 \item the dual equality: Does $\inf(P)=\sup(D)$ hold?
 \item the primal attainment: Does there exist a solution $\lh$ to
    $(P)?$ What about the minimizing sequences?
 \item the dual attainment: Does there exist a solution $\yb$ to $(D)?$
 \item the representation of the primal solutions: Find an identity
    of the type: $\lh\in\partial\F(T^\ast \yb).$
\end{itemize}
\end{questions}
In the case where the constraint set $C=\{x\}$ is reduced to a
single point,  the value $\sup (D)$ of the dual problem is
\begin{equation*}
    \Ls(x):=\sup_{y\in\YY}\{\yx - \F(T^\ast y)\}, \quad x\in\XX
\end{equation*}
which is the convex conjugate of $\La(y):=\F(T^\ast y), y\in\YY.$
\\
 We are going to answer the above questions in terms of
some \emph{extension} $\Fb$ of $\F$ under the weak constraint
qualification
\begin{equation}\label{eq-14b}
    C\cap\diffdom\Ls\not =\emptyset
\end{equation}
where
    $
\diffdom\Ls=\{x\in\XX;\partial_{\XX^\ast}\Ls(x)\not=\emptyset\}
    $
is the subset of all vectors in  $\XX$ at which $\Ls$ admits a
nonempty subdifferential with respect to the algebraic dual
pairing $\langle\XX,\XX^\ast\rangle$ with $\XX^\ast$ the algebraic
dual space of $\XX.$ Note that by the geometric version of
Hahn-Banach theorem, the intrinsic core of  $\Ls:$ $\icordom\Ls,$
is included in $\diffdom\Ls.$ Hence, a useful criterion to get
(\ref{eq-14b}) is
\begin{equation}\label{eq-14}
  C\cap\icordom\Ls\not =\emptyset.
\end{equation}
The drawback of such a general approach is that one has to compute
the extension $\Fb.$ In specific examples, this might be a
difficult task. In the case of the Monge-Kantorovich problem
\cite{p-Leo07b} it is immediate, but it requires some work in the
case of entropy minimization \cite{p-Leo07c}.

The restriction (\ref{eq-14}) is very weak since the intrinsic
core is the notion of interior which gives the largest possible
set. As $C\cap\dom\Ls=\emptyset$ implies that $(P)$ has no
solution, the only case where the problem remains open when
$\icordom\Ls$ is nonempty is the situation where $C$ and $\dom\Ls$
are tangent to each other. This is used in \cite{p-Leo07c} to
obtain general results for convex integral functionals. The
representation of their minimizers, see (\ref{eq-22}), are
obtained under the constraint qualification (\ref{eq-14}) which is
much weaker than the usual constraint qualification: $$\inter
C\cap\dom\Lambda^*\not=\emptyset$$ where $\inter C$ is the
interior of $C$ with respect to some topology which is not
directly connected to the ``geometry'' of $\Lambda^*.$ In
particular, $\inter C$ must be nonempty; this is a considerable
restriction.
\\
Nevertheless, the Monge-Kantorovich optimal transport problem
provides an interesting case where the specifications of the
constraints never stand in $\icordom\Ls,$ see \cite[Remark
\remMK]{p-Leo07b}, so that (\ref{eq-14}) is useless and
(\ref{eq-14b}) is the right assumption to be used.

\subsection*{The strategy} A usual way to prove the dual
attainment and obtain some representation of the primal solutions
is to require that the constraint is qualified: a property which
allows to separate the convex constraint set $T^{-1}(C)$ and the
level sets of the objective function $\Fs.$ The strategy of this
article is different: one chooses suitable topologies so that the
level sets have nonempty interiors. This also allows to apply
Hahn-Banach theorem, but this time the constraint set is not
required to have a nonempty interior. We take the rule not to
introduce arbitrary topological assumptions since $(P)$ is
expressed without any topological notion. Because of the convexity
of the problem, one takes advantage of geometric easy properties:
the topologies to be considered later are associated with
seminorms which are gauges of level sets of the convex functions
$\F$ and $\Fs.$ They are useful tools to work with the
\emph{geometry} of $(P).$
\\
It appears that when the constraints are infinite-dimensional one
can choose several different spaces $\YY$ without modifying the
value and the solutions of $(P).$ Consequently, for a small space
$\YY$ the dual attainment is not the rule. As a consequence, we
are facing the problem of finding an \emph{extension} of $(D)$
which admits solutions in generic cases and such that the
representation of the primal solution is $\lh\in\partial\Fb(T^\ast
\yb)$ where $\Fb$ is some extension of $\F.$
\\
We are going to
\begin{itemize}
    \item use the standard saddle-point approach to convex problems
    based on conjugate duality as developed by Rockafellar in \cite{Roc74}
    \item with topologies which reflect some of the geometric
    structure of the objective function.
\end{itemize}
These made-to-measure topologies are associated with the gauges of
the level sets of $\F$ and $\Fs.$

\subsection*{Outline of the paper}
The results are stated without proof at Section
\ref{sec:abstractpb}. Their proofs are postponed to Section
\ref{sec:proofs}. Examples are introduced at Section
\ref{sec:examples} where one considers the Monge-Kantorovich
transport and entropy minimization problems. These problems are
investigated in \cite{p-Leo07b, p-Leo07c}.

\subsection*{Notation}
Let $X$ and $Y$ be topological vector spaces. The algebraic dual
space of $X$ is $X^{\ast},$ the topological dual space of $X$ is
$X'.$ The topology of $X$ weakened by $Y$ is $\sigma(X,Y)$ and one
writes $\langle X,Y\rangle$ to specify that $X$ and $Y$ are in
separating duality.
\\
Let $f: X\rightarrow [-\infty,+\infty]$ be an extended numerical
function. Its convex conjugate with respect to $\langle
X,Y\rangle$ is $f^*(y)=\sup_{x\in X}\{\langle x,y\rangle
-f(x)\}\in [-\infty,+\infty],$ $y\in Y.$ Its subdifferential at
$x$ with respect to $\langle X,Y\rangle$ is $\partial_Y
f(x)=\{y\in Y; f(x+\xi)\geq f(x)+\langle y,\xi\rangle, \forall
 \xi\in X\}.$ If no confusion occurs, one writes
$\partial f(x).$
\\
The intrinsic core of a subset $A$ of a vector space is
    $
\mathrm{icor\,} A=\{x\in A; \forall x'\in\mathrm{aff\,}A, \exists
t>0, [x,x+t(x'-x)[\subset A\}
    $
where $\mathrm{aff\,}A$ is the affine space spanned by $A.$
$\icordom f$ is the intrisic core of the effective domain of $f$
$\dom f=\{x\in X; f(x)<\infty\}.$
\\
The indicator of a subset $A$ of $X$ is defined by
    \begin{equation}\label{eq-15}
    \iota_A(x)=\left\{%
\begin{array}{ll}
    0, & \hbox{if }x\in A \\
    +\infty, & \hbox{otherwise} \\
\end{array}%
\right.,\quad x\in X.
    \end{equation}
The support function of $A\subset X$ is $\iota_A^*(y)=\sup_{x\in
A}\langle x,y\rangle,$ $y\in Y.$

\section{Statements of the results}\label{sec:abstractpb}

The dual equality and the primal attainment are stated at Theorem
\ref{xP3}; the dual attainment and the dual representation of the
minimizers are stated at Theorems \ref{T3a} and \ref{T3b}. Their
proofs are postponed to Section \ref{sec:proofs}.

\subsection{Basic diagram}\label{sec:Tadj}

Let $\UUo$ be a vector space, $\LLo=\UUo^\ast$ its algebraic dual
space, $\Fo$ a $(-\infty,+\infty]$-valued convex function on
$\UUo$ and $\Fos$ its convex conjugate for the duality $\langle
\UUo,\LLo\rangle:$
\begin{equation*}
    \Fos(\ell):= \sup_{u\in\UUo}\{\ul-\Fo(u)\},\quad \ell\in\LLo\\
\end{equation*}
 Let $\YYo$ be another vector space,
$\XXo=\YYo^\ast$ its algebraic dual space and
$\To:\LLo\rightarrow\XXo$ is a linear operator. We consider the
convex minimization problem
\begin{equation}
  \textsl{minimize } \Fos(\ell) \textsl{ subject to } \To\ell\in \Co,\
  \ell\in\LLo
\tag{$P_o$}
\end{equation}
where $\Co$ is a convex subset of $\XXo.$
\\
It is useful to define the constraint operator $\To$ by means of
its adjoint $\To^\ast :\YYo\rightarrow\LLo^\ast$ as follows. For
all $\ell\in\LLo,$
\begin{displaymath}
  \langle y,\To\ell\rangle_{\YYo,\XXo}=\langle\To^\ast y,\ell\rangle_{\LLo^\ast,\LLo},\quad\forall y\in\YYo.
\end{displaymath}
We shall assume that the restriction
\begin{equation}\label{eq-01}
    \To^\ast (\YYo)\subset \UUo
\end{equation}
 holds, where $\UUo$ is identified with a subspace of
$\LLo^\ast=\UUo^{\ast\ast}.$ It follows that the diagram
\begin{equation}
 \begin{array}{ccc}
\Big\langle\ \UUo & , & \LLo \ \Big\rangle \\
\To^\ast  \Big\uparrow & & \Big\downarrow
 \To
\\
\Big\langle\ \YYo & , & \XXo\ \Big\rangle
\end{array}
\tag{Diagram 0}
\end{equation}
is meaningful.

\subsection{Assumptions}
\label{sec:assumptions}

Let us give the list of our main hypotheses.
\begin{itemize}
\item[\HF]
1-\quad $\Fo: \UUo\rightarrow [0,+\infty]$ is $\sUL$-\lsc, convex and $\Fo(0)=0$\\
2-\quad $\forall u\in\UUo, \exists \alpha >0, \Fo(\alpha u)<\infty$\\
3-\quad $\forall u\in\UUo, u\not=0,\exists t\in\R, \Fo(tu)>0$
\item[\HT]
1-\quad $\To^\ast (\YYo)\subset\UUo$\\
2-\quad $\mathrm{ker\ }\To^\ast =\{0\}$ \item [\HC] \qquad $\Ci:=
 \Co\cap\XX$ is a convex $\sigma(\XX,\YY)$-closed subset of $\XX$
\end{itemize}
The definitions of the vector spaces $\XX$ and $\YY$ which appear
in the last assumption are stated below at Section
\ref{sec:problems}. For the moment, let us only say that if $ \Co$
is convex and $\sigma(\XXo,\YYo)$-closed, then \HC\ holds.

\par\smallskip
\noindent\textit{Comments about the assumptions.}\
\begin{itemize}
    \item[-] By construction, $\Fos$ is a convex $\sigma(\LLo,\UUo)$-closed
function, even if $\Fo$ doesn't satisfy \HFi. Assuming this
property of $\Fo$ is not a restriction.
    \item[-] The assumption \HFi\ also expresses that $\Fo$ achieves its minimum
at $u=0$ and that $\Fo(0)=0.$ This is a practical normalization
requirement which will allow us to build a gauge functional
associated with $\Fo.$ More, \HFi\ implies that $\Fos$ also shares
this property. Gauge functionals related to $\Fos$ will also
appear later.
    \item[-] With any convex function $\tilde\F$ satisfying \HFii, one can
associate a function $\Fo$ satisfying \HFi\ in the following
manner. Because of \HFii, $\tilde\F(0)$ is finite and there exists
$\ell_o\in\LLo$ such that $\ell_o\in\partial\tilde\F(0).$ Then,
the closed convex regularization $\Fo$ of
$u\in\UUo\mapsto\tilde\F(u)-\langle\ell_o,u\rangle -\tilde\F(0),$
satisfies \HFi\ and
${\tilde\F}^*(\ell)=\Fos(\ell-\ell_o)-\tilde\F(0),$ $\ell\in\LLo.$
    \item[-] The hypothesis \HFiii\ is not a restriction. Indeed, assuming \HFi, let us
suppose that there exists a direction $u_o\not=0$ such that
$\Fo(tu_o)=0$ for all real $t.$ Then any $\ell\in\LLo$ such that
$\langle\ell,u_o\rangle\not=0$ satisfies
$\Fos(\ell)\geq\sup_{t\in\R}t\langle\ell,u_o\rangle=+\infty$ and
can't be a solution to $(P).$
    \item[-] The hypothesis \HTii\ isn't a restriction either: If
$y_1-y_2\in\mathrm{ker\ }\To^\ast ,$ we have $\langle
\To\ell,y_1\rangle=\langle \To\ell,y_2\rangle,$ for all
$\ell\in\LLo.$ In other words, the spaces $\YYo$ and
$\YYo/\mathrm{ker\ }\To^\ast $ both specify the same constraint
sets $\{\ell\in\LLo; \To\ell=x\}.$
\end{itemize}
The effective assumptions are the following ones.
\begin{itemize}
    \item[-] The specific form of the objective function $\Fos$ as a
convex conjugate makes it a convex $\sLU$-closed function.
    \item[-] \HFii\ and \HC\ are geometric restrictions.
    \item[-] \HTi\ is a
regularity assumption on $\To.$
\end{itemize}

\subsection{Variants of $(P)$ and $(D)$}
\label{sec:problems} These variants are expressed below in terms
of new spaces and functions. Let us first introduce them.

\subsubsection*{The norms $|\cdot|_\F$ and $|\cdot|_\La$}
Let $\F_{\pm}(u)=\max(\Fo(u),\Fo(-u)).$ By \HFi\ and \HFii,
$\{u\in\UUo; \F_{\pm}(u)\leq 1\}$  is a convex absorbing balanced
set. Hence its gauge functional which  is defined for all $
u\in\UUo$ by $|u|_\F := \inf\{\alpha>0; \F_{\pm}(u/\alpha))\leq
1\}$ is a seminorm. Thanks to hypothesis \HFiii, it is a norm.
\\
Taking \HTi\ into account, one can define
\begin{equation}
\label{La}
  \Lao(y):= \Fo(\To^\ast y), y\in\YYo.
\end{equation}
Let $\La_\pm(y)= \max(\Lao(y),\Lao(-y)).$ The gauge functional on
$\YYo$ of the set $\{y\in\YYo; \La_\pm(y)\leq 1\}$ is $|y|_\La :=
\inf\{\alpha>0;\La_\pm(y/\alpha)\leq 1\}, y\in\YYo.$ Thanks to
\HF\ and \HT, it is a norm and
\begin{equation}\label{eq-117}
    |y|_\La=|\To^\ast y|_\F, \quad y\in\YYo.
\end{equation}

\subsubsection*{The spaces}
Let
\begin{eqnarray*}
& &\UU \mathrm{\ be\ the\ } \NF\textrm{-completion\ of\ } \UUo\mathrm{\ and\ let }\\
  & &\LL:= (\UUo,\NF)'\textrm{\ be\ the\ topological\ dual\ space\ of\ }
  (\UUo,\NF).
\end{eqnarray*}
Of course, we have $(\UU,\NF)'\cong\LL\subset\LLo$ where any
$\ell$
 in $\LL$ is identified with its restriction to $\UUo.$ Similarly, we introduce
\begin{eqnarray*}
& &\YY \mathrm{\ the\ } \NL\textrm{-completion\ of\ } \YYo
  \mathrm{\ and }\\
 & &\XX:= (\YYo,\NL)'\mathrm{\ the\ topological\ dual\ space\ of\
 }(\YYo,\NL).
\end{eqnarray*}
We have $(\YY,\NL)'\cong\XX\subset\XXo$ where any $x$ in $\YY'$ is
identified with its restriction to $\YYo.$
\\
We also have to consider the \emph{algebraic} dual spaces
$\LL^\ast$ and $\XX^\ast$ of $\LL$ and $\XX$.

\subsubsection*{The operators $T$ and $T^\ast$} It will be proved at Lemma
\ref{L2} that
\begin{equation}\label{eq-101}
    \To\LL\subset\XX
\end{equation}
Let us denote $T$ the restriction of $\To$ to $\LL\subset\LLo.$ By
(\ref{eq-101}), we have $T: \LL\to\XX.$ Let us define its adjoint
$T^\ast :\XX^\ast\rightarrow\LL^\ast$ for all $\omega\in\XX^\ast$
by:
\begin{displaymath}
  \langle  \ell,T^\ast\omega\rangle_{\LL,\LL^\ast}=\langle T\ell,\omega\rangle_{\XX,\XX^\ast},
  \forall \ell\in\LL.
\end{displaymath}
This definition is meaningful, thanks to (\ref{eq-101}). It will
be proved at Lemma \ref{L2} that
\begin{equation}\label{eq-102}
    T^*\YY\subset\UU
\end{equation}
We have the inclusions $\YYo\subset\YY\subset\XX^\ast.$ The
adjoint operator $\To^*$ is the restriction of $T^*$ to $\YYo.$

\subsubsection*{Some modifications of $\Fo$ and $\Lao$}

We introduce the following modifications of $\Fo:$
\begin{eqnarray*}
    \F(u)&:=&\sup_{\ell\in\LL}\{\ul-\Fos(\ell)\},\quad u\in\UU\\
     \Fb(\zeta)&:=&\sup_{\ell\in\LL}\{\lz-\Fs(\ell)\},\quad \zeta\in\LL^\ast.
\end{eqnarray*}
They are respectively $\sULi$ and $\sULii$-closed convex
functions. It is immediate to see that the restriction of $\Fb$ to
$\UU$ is $\F.$ As $\LL=\UU',$ $\F$ is also the
$|\cdot|_{\F}$-closed convex regularization of $\Fo.$ The function
$\Fb$ is the extension which appears in the introductory Section
\ref{sec:introduction}. We also introduce
\begin{eqnarray*}
  \La(y)&:=& \F(T^\ast y),\quad  y\in\YY\\
   \Lb(\omega)&:=&\Fb(T^\ast \omega),\quad \omega\in\XX^\ast
\end{eqnarray*}
which look like the definition (\ref{La}). Note that thanks to
(\ref{eq-102}), the first equality is meaningful. Because of the
previous remarks, the restriction of $\Lb$ to $\YY$ is $\La.$

\subsubsection*{The optimization problems}
Let  $\Fos$ and $\F^*$ be the convex conjugates of $\Fo$ and $\F$
with respect to the dual pairings $\langle\UUo,\LLo\rangle$ and
$\langle\UU,\LL\rangle:$
\begin{eqnarray*}
    \Fos(\ell)&:=& \sup_{u\in\UUo}\{\ul-\Fo(u)\},\quad \ell\in\LLo\\
    \Fs(\ell)&:=& \sup_{u\in\UU}\{\ul-\F(u)\},\quad \ell\in\LL
\end{eqnarray*}
and $\Los, \Ls$ be the convex conjugates of  $\Lao,\La$ with
respect to the dual pairings $\langle\YYo,\XXo\rangle$ and
$\langle\YY,\XX\rangle:$
\begin{eqnarray*}
     \Los(x)&:=& \sup_{y\in\YYo}\{\yx-\Lao(y)\},\quad x\in\XXo\\
    \Ls(x)&:=& \sup_{y\in\YY}\{\yx-\La(y)\},\quad x\in\XX\\
\end{eqnarray*}
 Finally,  denote
    $
  \Ci= \Co\cap\XX.
    $\\
 The optimization problems to be considered are:
\begin{align}
&\textsl{minimize } \Fos(\ell) & &\textsl{subject to } \To\ell\in
\Co,&& \ell\in\LLo   \tag{$P_o$}\\
&\textsl{minimize } \Fs(\ell) & &\textsl{subject to }
T\ell\in\Ci,&& \ell\in\LL   \tag{$P$}\\
&\textsl{minimize } \Ls(x) & &\textsl{subject to }
x\in\Ci,&& x\in\XX   \tag{$P^\XX$}\\
&\textsl{maximize } \inf_{x\in \Co}\yx - \Lao(y), && &
&y\in\YYo  & &   \tag{$D_o$}\\
&\textsl{maximize } \inf_{x\in\Ci}\yx - \La(y), && &
&y\in\YY  & &   \tag{$D$}\\
&\textsl{maximize } \inf_{x\in\Ci}\xo - \Lb(\omega), && &
&\omega\in\XX^\ast & & \tag{$\overline{D}$}
\end{align}

\subsection{Statements}

We are now ready to give answers to the Questions \ref{questions}
related to $(P)$ and $(D).$

\begin{theorem}[Primal attainment and dual equality]\label{xP3}
Assume that \HF\ and \HT\ hold.
\begin{enumerate}
    \item[(a)] For all $x$ in $\XXo,$ we have the little dual equality
        \begin{equation}\label{xped}
            \inf\{\Fos(\ell); \ell\in\LLo, \To\ell=x\}=\Los(x)\in
            [0,\infty].
        \end{equation}
        Moreover, in restriction to $\XX,$ $\Los=\Ls$ and
        $\Ls$ is $\sigma(\XX,\YY)$-inf-compact.
    \item[(b)] The problems \Po\ and \Pi\ are equivalent: they have
        the same solutions and $\inf(P_o)=\inf(P)\in[0,\infty].$
    \item[(c)] If $ \Co$ is convex and $\sXY$-closed, we have the dual equality
\begin{displaymath}
  \inf(P_o)=\sup(D_o)\in[0,\infty].
\end{displaymath}
\end{enumerate}
Assume that \HF, \HT\ and \HC\ hold.
\begin{enumerate}
    \item[(d)]  We have the dual equalities
\begin{equation}
   \inf(P_o)=\inf(P)=\sup(D)=\sup(\overline{D})=\inf_{x\in \Co}\Los(x)=\inf_{x\in \Ci}\Ls(x)\in [0,\infty]\label{xed1}
\end{equation}
    \item[(e)] If in addition $\inf(P_o)<\infty,$ then \Po\ is attained in
    $\LL.$ Moreover, any minimizing sequence of \Po\ has $\sLUi$-cluster
    points and every such cluster point solves \Po.
    \item[(f)] Let $\lh\in\LL$ be a solution to $(P),$ then
$\xh:= \To\lh$ is a solution to \PX\ and
$\inf(P)=\Fs(\lh)=\Los(\xh).$
\end{enumerate}
\end{theorem}

\begin{theorem}[Dual attainment and representation. Interior convex
constraint]\label{T3a}
    Assume that \HF, \HT\ and \HC\ hold.
\begin{enumerate}[(1)]
    \item
For any  $\lh\in\LL$ and $\ob\in\XX^\ast,$
 \begin{equation}\label{xeq-96}
    \left\{\begin{array}{cl}
      (a) & T\lh\in \Ci \\
      (b) & \langle \ob,T\lh\rangle\leq \langle \ob,x\rangle
      \textrm{ for all }x\in \Ci \\
      (c) & \lh\in\partial_{\LL}\Fb(T^*\ob) \\
    \end{array}\right.
\end{equation}
is equivalent to
\begin{equation*}
    \left\{%
\begin{array}{l}
    \hbox{$\lh$ is a solution to \Po\ and} \\
    \hbox{$\ob$ is a solution to \Dii} \\
\end{array}%
\right.
\end{equation*}

    \item Suppose that in addition the interior constraint qualification
 \begin{equation}\label{xeq-106}
    \Co\cap \mathrm{icor\,}(\To\dom\Fos)\not=\emptyset
\end{equation}
is satisfied. Then, the primal problem \Po\ is attained in $\LL$
and the dual problem \Dii\ is attained in
    $\XX^*.$
\end{enumerate}
\end{theorem}
Note that (\ref{xeq-106}) is equivalent to $
\Co\cap\icordom\Los\not=\emptyset.$

As can be seen in \cite[Remark \remMK]{p-Leo07b}, the
Monge-Kantorovich problem provides an example where no constraint
is interior. In order to solve it, we are going to consider the
more general situation (\ref{eq-14b}) where the constraint is said
to be a \emph{subgradient constraint}. This means that $ \xh$
belongs to
\begin{eqnarray*}
    \diffdom\Los&=&\{x\in\XX; \partial_{\XX^\ast}\Los(x)\not
    =\emptyset\}\quad \textrm{ where}\\
    \partial_{\XX^\ast}\Los(x)&=&\{\omega\in\XX^\ast; \Los(x')\geq\Los(x)+\langle x'-x,\omega\rangle,\forall
    x'\in\XX\}.
\end{eqnarray*}
 Two new optimization problems to be considered are
\begin{align}
&\textsl{minimize } \Fos(\ell) & &\textsl{subject to } \To\ell=\xh,&& \ell\in\LLo   \tag{$P^{\xh}$}\\
&\textsl{maximize } \langle\xh,\omega\rangle-\Lb(\omega), & &&
&\omega\in\XX^\ast & & \tag{$\overline{D}^{\xh}$}
\end{align}
where $\xh\in\XXo.$ This corresponds to the simplified case where
$\Co$ is reduced to the single point $\xh.$

\begin{theorem}[Dual attainment and representation. Subgradient affine constraint]\label{T3b}
Let us assume that \HF\ and \HT\ hold.
\begin{enumerate}[(1)]
    \item
For any  $\lh\in\LL$ and $\ob\in\XX^\ast,$
 \begin{equation}\label{eq-97}
    \left\{\begin{array}{cl}
      (a) & T\lh=\xh\\
      (b) & \lh\in\partial_{\LL}\Fb(T^*\ob) \\
    \end{array}\right.
\end{equation}
is equivalent to
\begin{equation*}
    \left\{%
\begin{array}{l}
    \hbox{$\lh$ is a solution to $(P^{\xh})$ and} \\
    \hbox{$\ob$ is a solution to  \Dxh} \\
\end{array}%
\right.
\end{equation*}

    \item Suppose that in addition the subgradient constraint qualification
\begin{equation}\label{xCQbis}
    \xh\in\diffdom\Los,
\end{equation}
is satisfied. Then, the primal problem $(P^{\xh})$ is attained in
$\LL,$ and the dual problem \Dxh\ is attained in
    $\XX^*.$
\end{enumerate}
\end{theorem}

It is well-known that the representation formula (\ref{xeq-96}-c)
or (\ref{eq-97}-b):
        \begin{equation}\label{xeq-110}
    \lh\in\partial_{\LL}\Fb(T^*\ob)
\end{equation}
 is equivalent to
 \begin{equation*}
    T^*\ob\in\partial_{\LL^*}\Fs(\lh)
\end{equation*}
 and also equivalent to Young's identity
\begin{equation}\label{xeq-107}
   \Fs(\lh)+\Fb(T^*\ob)=\langle \ob,T\lh\rangle=\Ls(\xh)+\Lb(\ob).
\end{equation}
Formula (\ref{xeq-110}) can be made a little more precise by means
of the following regularity result.

\begin{theorem}
Assume that \HF, \HT\ and \HC\ hold. Any solution $\ob$ of \Dii\
or \Dxh\ shares the following properties
\begin{itemize}
 \item[(a)] $\ob$ is in the $\sigma(\XX^*,\XX)$-closure of $\dom\La;$
 \item[(b)] $T^\ast \ob$ is in the $\sigma(\LL^*,\LL)$-closure of
 $T^\ast(\dom\La).$
\end{itemize}
If in addition the level sets of $\F$ are $\NF$-bounded, then
\begin{itemize}
 \item[(a')] $\ob$ is in $\YY''.$ More precisely, it is in the $\sigma(\YY'',\XX)$-closure of $\dom\La;$
 \item[(b')] $T^\ast \ob$ is in $\UU''.$ More precisely, it is in the $\sigma(\UU'',\LL)$-closure of
 $T^\ast(\dom\La)$
\end{itemize}
where $\YY''$ and $\UU''$ are the \emph{topological} bidual spaces
of $\YY$ and $\UU.$ This occurs if $\F,$ and therefore $\Fs,$ is
an even function.
\end{theorem}

\section{Examples}\label{sec:examples}

The abstract results of Section \ref{sec:abstractpb} are
exemplified by means of the Monge-Kantorovich optimal transport
problem and the problem of minimizing entropy functionals on
convex sets.

\subsection{The Monge-Kantorovich optimal transport problem}

Denote $\PA, \PB$ and  $\PAB$ the sets of all probability measures
on the spaces $A,$ $B$ and $\AB.$ Let $c:\AB \to [0,\infty)$ a
nonnegative (cost) function and two probability measures
$\mu\in\PA$ and $\nu\in\PB$ on $A$ and $B.$ The Monge-Kantorovich
problem is
\begin{equation}\label{MK}
    \textsl{minimize } \pi\in\PAB\mapsto \IAB c(a,b)\,\pi(dadb)
    \textsl{ subject to }\pi\in P(\mu,\nu) \tag{MK}
\end{equation}
where $P(\mu,\nu)$ is the set of all $\pi\in\PAB$ with prescribed
marginals $\pi_A=\mu$ on $A$ and $\pi_B=\nu$ on $B.$ Any solution
of (MK) is called an optimal plan. For a general account on this
active field of research, see C.~Villani's book \cite{Vill03}.
\\
Without going into the details, let us indicate how this problem
enters the present framework. Denote $C_A,$ $C_B$ and $\CAB$ the
spaces of all continuous bounded functions on $A,$ $B$ and $\AB.$
The function $\Fo$ is  defined on the space $\UUo=\CAB$ by
$$
\Fo(u)=\iota_{\{u\le c\}},\quad u\in\CAB
$$
see (\ref{eq-15}). The marginal constraint $\pi\in P(\mu,\nu)$ is
obtained choosing $\YYo=C_A\times C_B$ and
$$
\To^\ast(f,g)=f\oplus g,\quad f\in C_A, g\in C_B
$$
with
$$
f\oplus g(a,b):=f(a)+g(b),\quad a\in A, b\in B,
$$
see Section \ref{sec:constraints} below. This gives
$\Lao(f,g)=\iota_{\{f\oplus g\le c\}}$ and the dual equality
(\ref{xped}) is the well-known Kantorovich dual equality
\begin{eqnarray*}
   && \inf\left\{ \IAB c(a,b)\,\pi(dadb);\pi\in P(\mu,\nu)\right\} \\
  &=& \sup\left\{\int_A f(a)\,\mu(da)+\int_B g(b)\,\nu(db); f\in C_A,
g\in C_B: f\oplus g\le c\right\}.
\end{eqnarray*}
In \cite{p-Leo07b}, cost functions $c$ which may take infinite
values are considered and Theorem \ref{T3b} is used to
characterize the optimal plans, yielding a new result on this
well-known optimization problem.

\subsection{Entropy minimization}\label{sec:entropmin}

The problem is sketched in this section and studied in further
details in \cite{p-Leo07c}.
\subsubsection*{Entropy}
Let $R$ be a positive measure on a space $\ZZ$ and take a
$[0,\infty]$-valued measurable function $\gamma^*$ on $\ZZ\times
\mathbb{R}$ such that $\gamma^*(z,\cdot):=\gamma_z^*$ is convex
and \lsc\ for all $z\in \ZZ.$ Denote $\MZ$ the space of all signed
measures $Q$ on $\ZZ.$ The entropy functional to be considered is
defined by
\begin{equation}\label{eq-16}
    I(Q)=
    \left\{
    \begin{array}{ll}
    \IZ \gamma^*_z(\frac{dQ}{dR}(z))\,R(dz) & \mathrm{if \ }
    Q\prec R\\
    +\infty & \mathrm{otherwise}
    \end{array}
\right.,\quad Q\in\MZ.
\end{equation}
where $Q\prec R$ means that $Q$ is absolutely continuous with
respect to $R.$ Assume that for each $z$ there exists a unique
$m(z)$ which minimizes $\gamma^*_z$ with $\gamma^*_z(m(z))=0.$
Then, $I$ is $[0,\infty]$-valued, its unique minimizer is $mR$ and
$I(mR)=0.$
\\
As for each $z,$ $\gamma^*_z$ is closed convex, it is the convex
conjugate of some closed convex function $\gamma_z.$ Defining
\begin{equation*}
    \lambda(z,s)=\gamma(z,s)-m(z)s,\quad z\in \ZZ, s\in\mathbb{R},
\end{equation*}
one sees that for $R$-a.e.\! $z,$ $\lambda_z$ is a nonnegative
convex function and it vanishes at 0. A favorable choice for
$\UUo$ is the space of all measurable functions $u$ on $\ZZ$ such
that
\begin{equation}\label{eq-50}
    \IZ \lambda(z,\alpha u(z))\,R(dz)<\infty, \textrm{ for all
    }\alpha\in \mathbb{R}.
\end{equation}
With
$$
\lmax(z,s)=\max[\lambda(z,s),\lambda(z,-s)]\in [0,\infty], \quad
z\in \ZZ, s\in \mathbb{R},
$$
(\ref{eq-50}) is equivalent to $u$ belongs to
    \[
    \UUo=\EEl:=\left\{u; \IZ \lmax(z,\alpha u(z))\,R(dz)<\infty,
\forall \alpha>0\right\}
    \]
the ``small" Orlicz space associated with the Young function
$\lmax.$  Taking
\begin{equation}\label{eq-18}
    \Fo(u)= \IZ\lambda_z(u(z))\,R(dz)\in [0,\infty],\quad u \in\EEl
\end{equation}
leads to
\begin{equation}\label{eq-20}
  I(Q)=\Fos(Q-mR),\quad Q\in\MZ
\end{equation}
This identity is a consequence of general results of Rockafellar
on conjugate duality for integral functionals \cite{Roc68}.
Moreover, the effective domain of $I$ is included in the space
\begin{equation*}
    \Ml:=\left\{Q\in\MZ; \IZ |u|\,d|Q|<\infty, \forall
    u\in\EEl\right\}.
\end{equation*}

\subsubsection*{Constraint}
In order to define the constraint, take $\XXo$ a vector space and
a function $\theta:\ZZ\rightarrow\XXo.$ One wants to give some
meaning to the formal constraint $ \IZ \theta\,dQ=x $ with
$Q\in\Ml$ and $x\in\XXo.$ Suppose that $\XXo$ is the algebraic
dual space of some vector space $\YYo$ and define for all
$y\in\YYo,$
\begin{equation}\label{eq-17}
  \To^\ast y(z):=\langle y,\theta(z)\rangle_{\YYo,\XXo},\
  z\in \ZZ.
\end{equation}
 Assuming that
\begin{equation*}
    \To^\ast y\in\EEl,\quad \forall y\in\YYo
\end{equation*}
allows to define the constraint operator
$$
\To Q:=\IZ\theta\,dQ,\quad Q\in\Ml
$$
by
\begin{equation*}
 \left\langle y, \IZ\theta \,dQ\right\rangle_{\YYo,\XXo} = \IZ \langle
y,\theta(z)\rangle_{\YYo,\XXo}\,Q(dz),\quad \forall y\in\YYo.
\end{equation*}

\subsubsection*{Minimization problem}
The entropy minimization problem to be considered is
\begin{equation*}
    \textsl{minimize } I(Q)
    \textsl{ subject to } \IZ\theta\,dQ\in \Co,\ Q\in\Ml
\end{equation*}
where $\Co$ is a convex subset of $\XXo.$

\subsubsection*{Results}

\newcommand{\Qh}{\widehat{Q}}

Applying the abstract results of the present paper, in
\cite{p-Leo07c} are obtained the following results. Let
$\Gamma^*(x)=\sup_{y\in\YYo}\{\yx-\IZ\gamma_z(\langle
y,\theta(z)\rangle)\,R(dz)\},$ $x\in\XXo.$ The dual equality is
$\inf\{I(Q); Q\in\Ml,  \IZ\theta\,dQ\in \Co\}=\inf_{\Co}\Gamma^* $
and under the assumption
\begin{equation}\label{eq-21}
    \Co\cap \icordom\Gamma^*\not=\emptyset,
\end{equation}
the characterization of the minimizer $\Qh$ is as follows.
Defining $\xh\eqdef \IZ\theta\,d\Qh$ in the weak sense with
respect to the duality $\langle\YYo,\XXo\rangle,$ $\Qh$ is a
minimizer if and only if there exists some linear form $\yb$ on
$\XXo$ such that $\langle\yb,\theta(\cdot)\rangle$ is measurable,
    $\IZ\lmax(z,\alpha_o\langle\yb,\theta(z)\rangle)\,R(dz)<\infty$  for some $\alpha_o>0$
and
 \begin{equation}\label{eq-22}
    \left\{\begin{array}{cl}
      (a) & \xh\in \Co\cap\dom\Gamma^* \\
      (b) & \langle \yb,\xh\rangle \leq \langle \yb,x\rangle, \forall x\in \Co\cap\dom\Gamma^* \\
      (c) & \Qh(dz)=\gamma'_z(\langle\yb, \theta(z)\rangle)\,R(dz). \\
    \end{array}\right.
\end{equation}
where $\gamma'_z(s)=\frac{\partial}{\partial s}\gamma(z,s).$

\begin{remark}\label{rem-01}
A usual form of constraint qualification required for this
representation is $\inter \Co\cap\dom\Gamma^*\not=\emptyset$ where
$\inter \Co$ is the interior of $\Co$ with respect to some
topology which is not directly connected to the ``geometry'' of
$\Gamma^*.$ In particular, $\inter \Co$ must be nonempty; this is
a considerable restriction. The constraint qualification $\Co\cap
\icordom\Gamma^*\not=\emptyset$ is much weaker.
\end{remark}

\subsubsection*{Literature about entropy minimization}
Entropy minimization problems appear in many areas of applied
mathematics and sciences. The literature about the minimization of
entropy functionals under convex constraints is considerable: many
papers are concerned with an engineering approach, working on the
implementation of numerical procedures in specific situations. In
fact, entropy minimization is a popular method to solve ill-posed
inverse problems.
\\
Surprisingly, rigorous general results on this topic are quite
recent. Let us cite, among others, the main contribution of
Borwein and Lewis: \cite{BL91a,BL91b,BL91c,BL92,BL93,BLN} together
with the paper \cite{TV} by Teboulle and Vajda. In these papers,
topological constraint qualifications are required: it is assumed
that the constraints stand in some topological interior of the
domain of $I.$ Such restrictions are removed in \cite{p-Leo07c}.

\subsection{Some examples of constraints}\label{sec:constraints}

Let us consider the two standard constraints which are the moment
constraints and the marginal constraints.

\subsubsection*{Moment constraints}
Let $ \theta=(\theta_k)_{1\leq k\leq K} $ be a measurable function
from $\ZZ$ to  $\XXo=\R^K.$ The moment constraint is defined by
$$
\IZ\theta\,dQ=\left(\IZ\theta_k\,dQ\right)_{1\leq k\leq K}\in\R^K,
$$
for each $Q\in\MZ$ which integrates all the real valued measurable
functions $\theta_k.$

\subsubsection*{Marginal constraints}
Let $\ZZ=\AB$ be a product space, $M_{AB}$ be the space of all
\emph{bounded} signed measures on $\AB$ and $U_{AB}$ be the space
of all measurable bounded functions $u$ on $\AB.$ Denote
$\ell_A=\ell(\cdot\times B)$ and $\ell_B=\ell(A\times\cdot)$ the
marginal measures of $\ell\in M_{AB}.$ The constraint of
prescribed marginal measures is specified by
$$
\IAB\theta\,d\ell=(\ell_A,\ell_B)\in M_{A}\times M_{B},\quad
\ell\in M_{AB}
$$
where $M_{A}$ and $M_{B}$ are the spaces of all bounded signed
measures on $A$ and $B.$ The function $\theta$ which gives the
marginal constraint is
$$
\theta(a,b)=(\delta_{a}, \delta_{b}),\ a\in A, b\in B
$$
where $\delta_a$ is the Dirac measure at $a.$ Indeed,
$(\ell_A,\ell_B)=\IAB (\delta_a,\delta_b)\,\ell(dadb).$
\\
More precisely, let $U_{A},$ $U_{B}$ be the spaces of measurable
functions on $A$ and $B$ and take $\YYo=U_{A}\times U_{B}$ and
$\XXo=U_{A}^*\times U_{B}^*.$ Then, $\theta$ is a measurable
function from $\ZZ=\AB$ to $\XXo=U_{A}^*\times U_{B}^*.$ It is
easy to see that the adjoint of the marginal operator
\begin{equation*}
    \To\ell=(\ell_A,\ell_B)\in U_{A}^*\times U_{B}^*,\quad
    \ell\in\LLo=U_{AB}^*
\end{equation*}
where $\langle f,\ell_A\rangle:=\langle f\otimes 1,\ell\rangle$
and $\langle g,\ell_B\rangle:=\langle 1\otimes g,\ell\rangle$ for
all $f\in U_A$ and $g\in U_B,$ is given by
\begin{equation}\label{eq-19}
    \To^\ast(f,g)=f\oplus g\in U_{AB},\quad f\in U_{A}, g\in U_{B}
\end{equation}
where $f\oplus g(a,b):=f(a)+g(b),$ $a\in A, b\in B.$

\section{Preliminary results}

In this section, one introduces notation and proves preliminary
technical results for the proofs of the results of Section
\ref{sec:abstractpb}.

\subsection{The saddle-point method (for fixing notation)}\label{sec:convexmin}

We are going to apply the general results of the Lagrangian
approach to the minimization problem \Po. To quote easily and
precisely some well-known results of convex minimization while
proving our results, we give a short overview of the approach to
convex minimization problems by means of conjugate duality as
developed in Rockafellar's monograph \cite{Roc74}.

Let $A$ be a vector space and $f: A\rightarrow [-\infty,+\infty]$
an extended real convex function. We consider the following convex
minimization problem
\begin{equation}
  \textsl{minimize } f(a), a\in A
\tag{$\mathcal{P}$}
\end{equation}
Let $Q$ be another vector space. The perturbation of the objective
function $f$ is a function $F:A\times Q\rightarrow
[-\infty,+\infty]$ such that for $ q=0\in Q,$ $F(\cdot,
0)=f(\cdot).$ The problem $(\mathcal{P})$ is imbedded in a
parametrized family of minimization problems
\begin{equation}
  \textsl{minimize } F(a,q), a\in A
\tag{$\mathcal{P}_q$}
\end{equation}
The value function of $(\mathcal{P}_q)_{q\in Q}$ is
\begin{displaymath}
  \varphi(q):= \inf(\mathcal{P}_q) =\inf_{a\in
  A}F(a,q)\in[-\infty,+\infty], q\in Q.
\end{displaymath}
Let us assume that the perturbation is chosen such that
\begin{equation}
  \label{Fcv}
  F \mathrm{\ is\ jointly\ convex\ on\ } A\times Q.
\end{equation}
Then, $(\mathcal{P}_q)_{q\in Q}$ is a family of convex
minimization problems and the value function $\varphi$ is convex.

Let $B$ be a vector space in dual pairing with $Q.$ This means
that $B$ and $Q$ are locally convex topological vector spaces in
separating duality such that their topological dual spaces $B'$
and $Q'$ satisfy $B'=Q$ and $Q'=B$ up to some isomorphisms. The
Lagrangian associated with the perturbation $F$ and the duality
$\langle B,Q\rangle$ is
\begin{equation}
\label{KF}
  K(a,b):= \inf_{q\in Q}\{\bq+F(a,q)\}, a\in A, b\in B.
\end{equation}
 Under (\ref{Fcv}), $K$ is a convex-concave function.
Assuming in addition that $F$ is chosen such that
\begin{equation}
  \label{Fcl}
 q\mapsto F(a,q) \mathrm{\ is\ a\ closed\ convex\ function\
for\ any\ } a\in A,
\end{equation}
one can reverse the conjugate duality relation (\ref{KF}) to
obtain
\begin{equation}
\label{FK}
  F(a,q)=\sup_{b\in B}\{K(a,b)-\bq\}, \forall a\in A, q\in Q
\end{equation}

Introducing another vector space $P$ in separating duality with
$A$ we define the function
\begin{equation}
\label{GK}
  G(b,p):= \inf_{a\in A}\{K(a,b)-\ap\}, b\in B, p\in P.
\end{equation}
This formula is analogous to (\ref{FK}). Going on symmetrically,
one interprets $G$ as the concave perturbation of the objective
concave function
\begin{displaymath}
  g(b):= G(b,0), b\in B
\end{displaymath}
associated with the concave maximization problem
\begin{equation}
  \textsl{maximize } g(b), b\in B
\tag{$\mathcal{D}$}
\end{equation}
which is  the dual problem of $(\mathcal{P}).$ It is imbedded in
the family of concave maximization problems $(\mathcal{D}_p)_{p\in
P}$
\begin{equation}
  \textsl{maximize } G(b,p), b\in B
\tag{$\mathcal{D}_p$}
\end{equation}
whose value function is
\begin{displaymath}
  \gamma(p):= \sup_{b\in B} G(b,p), p\in P.
\end{displaymath}
Since $G$  is jointly concave, $\gamma$ is also concave. We have
the following diagram
\begin{displaymath}
 \begin{array}{lrccclr}
       &           & \gamma(p)  &         & f(a)   &            &\\
       &\Big\langle& P          & ,       & A      &\Big\rangle &\\
G(b,p) &           &            & K(a,b)  &        &            & F(a,q)\\
       &\Big\langle& B          & ,       & Q      &\Big\rangle &\\
       &           & g(b)       &         & \varphi(q)&         &
\end{array}
\end{displaymath}

The concave conjugate of the function $f$ with respect to the dual
pairing $\langle Y,X\rangle$ is $f^\cav(y)=\inf_x\{\langle
y,x\rangle-f(x)\}$ and its superdifferential at $x$ is
$\widehat{\partial}f(x)=\{y\in Y; f(x')\leq f(x)+\langle y,
x'-x\rangle\}.$

\begin{theorem}
\label{TLag1} We assume that $\langle P,A\rangle$ and $\langle
B,Q\rangle$ are topological dual pairings.
\begin{itemize}
 \item[(a)] We have $\sup(\mathcal{D})=\varphi^{**}(0).$ Hence, the dual equality
$\inf(\mathcal{P})=\sup(\mathcal{D})$ holds if and only if
$\varphi(0)=\varphi^{**}(0).$
 \item[(b)] In
particular,
\[\left.
\begin{array}{l}
  \bullet\ F \textrm{ is jointly convex }  \\
  \bullet\ \varphi\textrm{ is lower semicontinuous at } 0 \\
  \bullet\ \sup(\mathcal{D}) >-\infty\\
\end{array}\right\}\Rightarrow \inf(\mathcal{P})=\sup(\mathcal{D})
\]
 \item[(c)] If the dual equality
holds, then
\begin{displaymath}
  \mathrm{argmax\ }g=-\partial \varphi(0).
\end{displaymath}
\end{itemize}

Let us assume in addition that (\ref{Fcv}) and (\ref{Fcl}) are
satisfied.
\begin{itemize}
 \item[(a')] We have $\inf(\mathcal{P})=\gamma^{\cav\cav}(0).$
Hence, the dual equality $\inf(\mathcal{P})=\sup(\mathcal{D})$
holds if and only if $\gamma(0)=\gamma^{\cav\cav}(0).$
 \item[(b')] In
particular,
\[\left.
\begin{array}{l}
  \bullet\ \gamma\textrm{ is \usc\ at } 0 \\
  \bullet\ \inf(\mathcal{P})<+\infty\\
\end{array}\right\}\Rightarrow \inf(\mathcal{P})=\sup(\mathcal{D})
\]
 \item[(c')] If the dual equality
holds, then
\begin{displaymath}
  \mathrm{argmin\ }f=-\widehat{\partial}\gamma(0).
\end{displaymath}
\end{itemize}
\end{theorem}

\begin{definition}[Saddle-point]
One says that $(\ab,\bb)\in A\times B$ is a \emph{saddle-point} of
the function $K$ if
\[
K(\ab,b)\leq K(\ab,\bb)\leq K(a,\bb),\quad\forall a\in A, b\in B.
\]
\end{definition}

\begin{theorem}[Saddle-point theorem and KKT relations]\label{TLag2}
The following statements are equivalent.

\begin{enumerate}[(1)]
    \item The point $(\ab,\bb)$ is a saddle-point of the Lagrangian $K$
    \item $f(\ab)\leq g(\bb)$
    \item The following  three statements hold
         \begin{enumerate}[(a)]
         \item we have the \emph{dual equality}: $\sup(\mathcal{D})= \inf(\mathcal{P}),$
         \item $\ab$ is a solution to the primal problem $(\mathcal{P})$ and
         \item $\bb$ is a solution to the dual problem $(\mathcal{D})$.
         \end{enumerate}
\end{enumerate}
In this situation, one also gets
\begin{equation}\label{eq-302}
    \sup(\mathcal{D})=
\inf(\mathcal{P})=K(\ab,\bb)=f(\ab)=g(\bb).
\end{equation}
Moreover, $(\ab,\bb)$ is a saddle-point of $K$ if and only if it
satisfies
\begin{eqnarray}
  \partial_a K(\ab,\bb)&\ni& 0 \label{KTa} \\
 \widehat{\partial}_b  K(\ab,\bb) &\ni& 0  \label{KTb}
\end{eqnarray}
where the subscript $a$ or $b$ indicates the unfixed variable.
\end{theorem}

\subsection{Gauge functionals associated with a convex
function}\label{sec:gauge}

The following result is probably well-known, but since I didn't
find a reference for it, I give its short proof.

Let $\theta: S\rightarrow[0,\infty]$ be an extended  nonnegative
convex function on a vector space $S,$ such that $\theta(0)=0.$
Let $\Ss$ be the algebraic dual space of $S$ and $\ts$ the convex
conjugate of $\theta:$
 \begin{displaymath}
   \ts(r):= \sup_{s\in S}\{\langle r,s \rangle-\theta(s)\},\quad r\in\Ss.
 \end{displaymath}
It is easy to show that $\ts: \Ss\rightarrow[0,\infty]$ and
$\ts(0)=0.$ We denote $C_\theta:=\{\theta\leq 1\}$ and
$C_{\theta^*}:=\{\theta^*\leq 1\}$ the unit level sets of $\theta$
and $\theta^*.$ The gauge functionals to be considered are
\begin{eqnarray*}
  \jt(s)&:= & \inf\{\alpha>0 ;  s\in \alpha C_\theta\}
=\inf\{\alpha>0 ; \theta(s/\alpha)\leq 1\}\in [0,\infty], s\in S.\\
\jts(r)&:= &\inf\{\alpha>0 ;  r\in\alpha C_{\theta^*}\}
=\inf\{\alpha>0 ; \theta^*(r/\alpha)\leq 1\}\in [0,\infty], r\in
\Ss.
\end{eqnarray*}
As 0 belongs to $\Ct$ and $\Cts,$ one easily proves that $\jt$ and
$\jts$ are positively homogeneous. Similarly, as $\Ct$ and $\Cts$
are convex sets, $\jt$ and $\jts$ are convex functions.
\begin{proposition}\label{Pgauge}
    For all $r\in\Ss,$ we have
\begin{displaymath}
  \frac{1}{2}\jts(r)\leq \ds(r):= \sup_{s\in\Ct}\langle r,s\rangle
\leq 2 \jts(r).
\end{displaymath}
We also have
\begin{displaymath}
  \mathrm{cone\ }\dom\ts=\dom\jts=\dom\ds
\end{displaymath}
where $\mathrm{cone\ }\dom\ts$ is the convex cone (with vertex
$0$) generated by $\dom\ts.$
\end{proposition}
\proof $\bullet$\ Let us first show that $\ds(r)\leq 2\jts(r)$ for
all $r\in\Ss.$ For all $s\in\Ct$ and $\alpha>\jts(r),$ $\langle
r,s \rangle =\langle r/\alpha,s\rangle \alpha\leq
[\theta(s)+\ts(r/\alpha)]\alpha\leq (1+1)\alpha.$ Then, optimize
both sides of this inequality.
\\
$\bullet$\ Let us show that $\jts(r)\leq 2\ds(r).$ If
$\ds(r)=\infty,$ there is nothing to prove. So, let us suppose
that $\ds(r)<\infty.$ As $0\in\Ct,$ we have $\ds(r)\geq 0.$
\\
\emph{First case: } $\ds(r)>0.$\ For all $s\in S$ and
$\epsilon>0,$ we have $s/[\jt(s)+\epsilon]\in\Ct.$ It follows that
$\langle r/\ds(r),s\rangle =\langle r,s/[\jt(s)+\epsilon]\rangle
\frac{\jt(s)+\epsilon}{\ds(r)}\leq \ds(r)
\frac{\jt(s)+\epsilon}{\ds(r)}=\jt(s)+\epsilon.$ Therefore,
$\langle r/\ds(r),s\rangle\leq \jt(s),$ for all $s\in S.$
\\
If $s$ doesn't belong to $\Ct,$ then $\jt(s)\leq\theta(s).$ This
follows from the the assumptions on $\theta:$ convex function such
that $\theta(0)=0=\min\theta$ and the positive homogeneity of
$\jt.$ Otherwise, if $s$ belongs to $\Ct,$ we have $\jt(s)\leq 1.$
Hence, $
  \langle r/\ds(r),s\rangle \leq\max(1,\theta(s)), \forall s\in S.
$ On the other hand, there exists $s_o\in S$ such that
$\ts(r/[2\ds(r)])\leq \langle r/[2\ds(r)],s_o\rangle
-\theta(s_o)+1/2.$ The last two inequalities provide us with
$\ts(r/[2\ds(r)])\leq
\frac{1}{2}\max(1,\theta(s_o))-\theta(s_o)+\frac{1}{2}\leq 1$
since $\theta(s_o)\geq 0.$ We have proved that $\jts(r)\leq
2\ds(r).$
\\
\emph{Second case: } $\ds(r)=0.$\ We have $\langle r,s\rangle\leq
0$ for all $s\in\Ct.$ As $\dom\theta$ is a subset of the cone
generated by $\Ct,$ we also have for all $t>0$ and
$s\in\dom\theta,$ $\langle tr,s\rangle\leq 0.$ Hence $\langle
tr,s\rangle -\theta(s)\leq 0$ for all $s\in S$ and $\ts(tr)\leq
0,$ for all $t\geq 0.$ As $\ts\geq 0,$ we have $\ts(tr)=0,$ for
all $t\geq 0.$ It follows that $\jts(r)=0.$ This completes the
proof of the equivalence of $\jts$ and $\ds.$
\\
$\bullet$\ Finally, this equivalence implies that
$\dom\jts=\dom\ds$ and as $\ts(0)=0$ we have $0\in\dom\ts$ which
implies that $\mathrm{cone\ }\dom\ts=\dom\jts.$
\endproof

\subsection{Preliminary technical results}
Recall that $|u|_\F=\inf\{\alpha>0 ; \F_{\pm}(u/\alpha)\leq 1\}$
with $\F_{\pm}(u)=\max(\Fo(u),\Fo(-u)).$ Its associated dual
uniform norm is
$$
|\ell|_\F^*:= \sup_{u, |u|_\F\leq 1}|\ul|,\quad \ell\in\LL
$$
The topological dual space of $(\LL, |\cdot|_\F^*)$ is denoted by
$\UU'':$ the  bidual space of $(\UU,|\cdot|_\F).$
\\
Similarly, recall that $|y|_\La=\inf\{\alpha>0 ;
\La_{\pm}(y/\alpha)\leq 1\}$ with
$\La_{\pm}(y)=\max(\Lao(y),\Lao(-y)).$ Its associated dual uniform
norm is
$$
|x|_\La^*:= \sup_{y, |y|_\La\leq 1}|\yx|, \quad x\in\XX
$$
The topological dual space of $(\XX, |\cdot|_\La^*)$ is denoted by
$\YY'':$ the  bidual space of $(\YY,|\cdot|_\La).$

\begin{lemma}\label{L2}
Let us assume \HF\ and \HT.
\begin{enumerate}[(a)]
    \item $\dom\Fos\subset\LL$ and $\dom\Los\subset\XX$
    \item $\To(\dom\Fos)\subset\dom\Los$ and $\To\LL\subset\XX$
    \item $\To$ is $\sLU$-$\sXY$-continuous
    \item $T^\ast: \XX^\ast\to\LL^\ast$ is $\sYXii$-$\sULii$-continuous
    \item $T: \LL\to\XX$ is  $\NF^*$-$\NL^*$-continuous
    \item $T^\ast\YY''\subset \UU''$ where $\YY''$ and $\UU''$ are
    the topological bidual spaces of $\YY$ and $\UU$
    \item $T^\ast \YY\subset\UU$ and $T^\ast: \YY\to\UU$ is $\sYXi$-$\sULi$-continuous
    \item $T: \LL\to\XX$ is $\sLUi$-$\sXYi$-continuous
\end{enumerate}
\end{lemma}

\proof \boulette{a} For all $\ell\in\LLo$ and $\alpha>0,$ Young's
inequality yields: $\ul= \alpha\langle\ell,u/\alpha\rangle\leq
[\Fo(u/\alpha)+\Fos(\ell)]\alpha,$ for all $u\in\UUo.$ Hence, for
any $\alpha>|u|_\F,$ $\ul\leq [1+\Fos(\ell)]\alpha.$ It follows
that $\ul\leq [1+\Fos(\ell)]|u|_\F.$ Considering $-u$ instead of
$u,$ one gets
\begin{equation}
  \label{eq-02}
  |\ul|\leq [1+\Fos(\ell)]|u|_\F, \forall u\in\UUo, \ell\in\LLo.
\end{equation}
It follows that $\dom\Fos\subset\LL.$ One proves
$\dom\Los\subset\XX$ similarly.

\Boulette{b}
    It is easy to show that
$
  \Los(\To\ell)\leq\Fos(\ell),
  $
for all $\ell\in\LLo.$ It follows immediately that
$\To(\dom\Fos)\subset\dom\Los.$
\\
 Let us consider
$|\cdot|_{\Fs_{\pm}}$ and $|\cdot|_{\Ls_{\pm}}$ the gauge
functionals of the level sets $\{\Fs_{\pm}\leq 1\}$ and
$\{\Ls_{\pm}\leq 1\}.$ As above,
\begin{equation}
  \label{eq-11}
  \Ls_{\pm}(\To\ell)\leq\Fs_{\pm}(\ell), \quad \forall \ell\in\LLo
\end{equation}
Therefore, $\To(\dom\Fs_{\pm})\subset\dom\Ls_{\pm}.$ On the other
hand, by Proposition \ref{Pgauge}, the linear space spanned by
$\dom\Fs_{\pm}$ is $\dom |\cdot|_{\Fs_{\pm}}$ and the linear space
spanned by $\dom\Ls_{\pm}$ is $\dom |\cdot|_{\Ls_{\pm}}.$ But,
$\dom |\cdot|_{\Fs_{\pm}}=\dom |\cdot|_\F^*=\LL$ and $\dom
|\cdot|_{\Ls_{\pm}}=\dom |\cdot|_\La^*=\XX$ by Proposition
\ref{Pgauge} again. Hence, $\To\LL\subset\XX.$

 \Boulette{c}
 To prove that $\To$ is continuous, one has to show that
for any $y\in \YYo,$ $\ell\in \LLo\mapsto \langle
y,\To\ell\rangle\in \mathbb{R}$ is continuous. We get
$\ell\mapsto\langle y,\To\ell\rangle=\langle \To^*y, \ell\rangle$
which is continuous since \HTi\ gives $\To^*y\in \UUo.$

\Boulette{d} It is a direct consequence of $\To\LL\subset\XX.$ See
the proof of (c).

\Boulette{e} We know by Proposition \ref{Pgauge} that
$|\cdot|_{\Fs_{\pm}}\sim\NF^*$ and $|\cdot|_{\Ls_{\pm}}\sim\NL^*$
are equivalent norms on $\LL$ and $\XX$ respectively. For all
$\ell\in\LL,$ $|T\ell|_\La^*\leq 2|T\ell|_{\Ls_{\pm}}
=2\inf\{\alpha>0 ;\Ls_{\pm}(T\ell/\alpha)\leq 1\} \leq
2\inf\{\alpha>0 ; \Fs_{\pm}(\ell/\alpha)\leq 1\}.$ This last
inequality follows from (\ref{eq-11}). Going on, we get
$|T\ell|_\La^*\leq 2|\ell|_{\Fs_{\pm}}\leq 4|\ell|_\F^*,$ which
proves that $T$ shares the desired continuity property with
$\|T\|\leq 4.$

 \Boulette{f}
 Let us take  $\omega\in\YY''.$ For all $\ell\in\LL,$
$|\langle T^\ast\omega,\ell\rangle_{\LL^\ast,\LL}|=|\langle
\omega,T\ell\rangle_{\YY'',\XX}|$ $\leq \|\omega\|_{\YY''}
|T\ell|_\La^*\leq \|\omega\|_{\YY''}\|T\| |\ell|_\F^*$ where
$\|T\|<\infty,$ thanks to (e). Hence, $T^\ast\omega$ stands in
$\UU''.$

 \Boulette{g}
Take  $y\in\YY.$  Let us show that $T^*y$ is the strong limit of a
sequence in $\UUo.$ Indeed, there exists a sequence $(y_n)$ in
$\YYo$ such that $\lim_{n\rightarrow\infty} y_n=y$ in $(\YY,\NL).$
Hence, for all $\ell\in\LL,$ $|\langle T^\ast y_n-T^\ast
y,\ell\rangle_{\LL^\ast,\LL}|=|\langle
y_n-y,T\ell\rangle_{\YY,\XX}|$ $\leq \|T\| |y_n-y|_\La
|\ell|_\F^*$ and $\sup_{\ell\in\LL, |\ell|_\F^*\leq 1}|\langle
T^\ast y_n-T^\ast y,\ell\rangle|\leq \|T\| |y_n-y|_\La$ tends to 0
as $n$ tends to infinity, where $T^\ast y_n$ belongs to $\UUo$ for
all $n\geq 1$ by \HTi. Consequently, $T^\ast y\in\UU.$
\\
The continuity statement now follows from (d).

 \Boulette{h}
By (b), $T$ maps $\LL$ into $\XX$ and because of (g): $T^\ast
\YY\subset\UU.$ Hence, for all $y\in\YY,$ $\ell\mapsto \langle
T\ell,y\rangle_{\XX,\YY}=\langle\ell,T^\ast y\rangle_{\LL,\UU}$ is
$\sLUi$-continuous.
 This completes the proof of Lemma \ref{L2}.
 \endproof

Recall that $\Fos,$  $\Lao^*$ and $\La^*$ are the convex
conjugates of $\Fo,$  $\Lao$ and $\La$ for the dual pairings
$\langle\UUo,\LLo\rangle,$ $\langle\YYo,\XXo\rangle$ and
$\langle\YY,\XX\rangle.$

\begin{lemma}\label{res-05}
Under the hypotheses \HF\ and \HT, we have
\begin{enumerate}[(a)]
    \item $\Fo=\F \textrm{ on } \UUo;$
    \item $\Lao=\La \textrm{ on }  \YYo ;$
    \item $\Fos=\F^* \textrm{ on }  \LL.$
\end{enumerate}
\end{lemma}
\proof (a) follows directly from Lemma \ref{L2}-a and the
assumption that $\Fo$ is closed convex. (b) follows from (a).
\\
Let us show (c). As $\UUo$ is a dense subspace of $\UU,$ we obtain
that $\F$ is the convex $\sigma(\UU,\LL)$-\lsc\ regularization of
$\Fo+\iota_{\UUo}$ where $\iota_{\UUo}$ is the convex indicator of
$\UUo.$ Since the convex conjugate of a function and the convex
conjugate of its convex \lsc\ regularization match, this implies
that $\Fos=\F^*$ on $\LL.$
\endproof

\begin{lemma}\label{L6}
Under the hypothesis \HF,
\begin{enumerate}[(a)]
    \item  $\Fos$ is $\sLU$-inf-compact and
    \item $\Fs$ is $\sLUi$-inf-compact.
\end{enumerate}

\end{lemma}
\proof \boulette{b}  Recall that we already obtained at
(\ref{eq-02}) that $
  |\ul|\leq [1+\Fos(\ell)]|u|_\F,
$ for all $u\in\UUo$ and $\ell\in\LLo.$ By completion, one deduces
that for all $\ell\in\LL$ and $u\in\UU,$ $|\ul|\leq
[1+\Fs(\ell)]|u|_\F$ ($\Fos=\Fs$ on $\LL,$ Lemma \ref{res-05}-c.)
Hence, $\Fs(\ell)\leq A$ implies that $|\ell|_\F^*\leq A+1.$
Therefore, the level set $\{\Fs\leq A\}$ is relatively
$\sLUi$-compact.
\\
By construction, $\Fs$ is $\sLUi$-\lsc. Hence, $\{\Fs\leq A\}$ is
$\sLUi$-closed and $\sLUi$-compact.

 \Boulette{a}
As $\Fos=\Fs$ on $\LL$ (Lemma \ref{res-05}-c),
$\dom\Fos\subset\LL$ (Lemma \ref{L2}-a) and $\UUo\subset\UU,$ it
follows from the $\sLUi$-inf-compactness of $\Fs$ that $\Fos$ is
$\sLU$-inf-compact.
\endproof

\section{Proofs of the results of Section \ref{sec:abstractpb}}\label{sec:proofs}

The results of Section \ref{sec:abstractpb} are a summing up of
Proposition \ref{P1}, Lemma \ref{C4} , Proposition \ref{P3},
Corollary \ref{res-09}, Lemma \ref{L6b}, Proposition
\ref{res-06a}, Corollary \ref{res-06b}, Proposition \ref{res-07}
and Proposition \ref{res-08}.

\subsection{A first dual equality}
In this  section we only consider the basic spaces
$\UUo,\LLo,\YYo$ and $\XXo.$ Let us begin applying Section
\ref{sec:convexmin} with $\langle P,A\rangle=\langle
\UUo,\LLo\rangle$ and
 $\langle B,Q\rangle=\langle \YYo,\XXo\rangle$ and the topologies are
 the weak topologies $\sLU,$ $\sUL,$ $\sXY$ and $\sYX.$ The
 function to be minimized is
 $f(\ell)=\Fos(\ell)+\iota_{ \Co}(\To\ell),$ $\ell\in\LLo.$  The perturbation $F$
of $f$ is Fenchel's one:
\begin{displaymath}
  F_0(\ell,x)=\Fos(\ell)+\iota_{ \Co}(\To\ell+x),\quad \ell\in\LLo, x\in\XXo.
\end{displaymath}
We assume \HTi: $\To^\ast \YYo\subset\UUo,$ so that the duality
diagram is
\begin{equation}
 \begin{array}{ccc}
\Big\langle\ \UUo & , & \LLo \ \Big\rangle \\
\To^\ast  \Big\uparrow & & \Big\downarrow
 \To
\\
\Big\langle\ \YYo & , & \XXo\ \Big\rangle
\end{array}
\tag{Diagram 0}
\end{equation}
The analogue of $F$ for the dual problem is
\begin{displaymath}
  G_0(y,u):= \inf_{\ell,x}\{\yx-\ul+F_0(\ell,x)\}
=\inf_{x\in \Co}\yx-\Fo(\To^\ast y+u),\quad y\in\YYo,u\in\UUo
\end{displaymath}
The corresponding value functions are
\begin{eqnarray*}
\varphi_0(x) &=& \inf\{\Fos(\ell); \ell\in\LLo: \To\ell\in
\Co-x\},\quad x\in\XXo
\\
 \gamma_0(u)&=&\sup_{y\in\YYo}\{\inf_{x\in \Co}\yx-\Fo(\To^\ast
 y+u)\},\quad  u\in\UUo.
\end{eqnarray*}
The primal and dual problems are \Po\ and \Do.

\begin{lemma}\label{res-04}
Assuming \HF\ and \HTi, if $\Co$ is a $\sXY$-closed convex set,
$F_0$ is jointly closed convex on $\LLo\times \XXo.$
\end{lemma}

\proof
 As $\To$ is linear continuous (Lemma \ref{L2}-c) and  $\Co$ is closed
convex, $\{(\ell,x); \To\ell+x\in \Co\}$ is closed convex in
$\LLo\times \XXo.$ As $\Fos$ is closed convex on $\LLo,$ its
epigraph is closed convex in $\LLo\times \mathbb{R}.$  It follows
that $\epi F_0=(\XXo\times\epi \Fos)\cap [\{(x,\ell); \To\ell+x\in
\Co\}\times\mathbb{R}]$ is closed convex, which implies that $F_0$
is convex and \lsc. As it is nowhere equal to $-\infty$ (since
$\inf F_0\geq \inf \Fos>-\infty$), $F_0$ is also a closed convex
function.
\endproof

Therefore, assuming that $ \Co$ is a $\sXY$-closed convex set, one
can apply the general theory of Section \ref{sec:convexmin} since
the perturbation function $F_0$ satisfies the assumptions
(\ref{Fcv}) and (\ref{Fcl}).

\begin{proposition}\label{P1}
Let us assume that \HF\ and \HT\ hold. If $ \Co$ is convex and
$\sXY$-closed, we have the dual equality
\begin{equation}\label{ed0}
  \inf(P_o)=\sup(D_o)\in[0,\infty].
\end{equation}
In particular,  for all $x$ in $\XXo,$ we have the little dual
equality
        \begin{equation}\label{ped0}
            \inf\{\Fos(\ell); \ell\in\LLo, \To\ell=x\}=\Los(x)\in
            [0,\infty].
        \end{equation}
\end{proposition}

\proof The identity (\ref{ped0}) is a special case of (\ref{ed0})
with $\Co=\{x\}.$
\\
To prove (\ref{ed0}), we consider separately the cases where
$\inf(P_o)<+\infty$ and $\inf(P_o)=+\infty.$

\par\medskip\noindent\textit{Case where $\inf(P_o)<+\infty.$}\
Thanks to Theorem \ref{TLag1}-b', it is enough to prove that
$\gamma_0$ is \usc\ at $u=0.$ We are going to prove that
$\gamma_0$ is continuous at $u=0.$ Indeed, for all $u\in\UUo,$
\[
-\gamma_0(u)=\inf_y\{\Fo(\To^*y+u)-\inf_{x\in \Co}\yx\}\leq \Fo(u)
\]
where the inequality is obtained taking $y=0.$ The norm $\NF$ is
designed so that $\Fo$ is  bounded above on a $\NF$-neighbourhood
of zero. By the previous inequality, so is the convex function
$-\gamma_0.$ Therefore, $-\gamma_0$ is $\NF$-continuous on
$\icordom (-\gamma_0)\ni 0.$ As it is convex and
$\LL=(\UUo,\NF)',$ it is also $\sigma(\UUo,\LL)$-\lsc\  and a
fortiori $\sUL$-\lsc\ , since $\LL\subset\LLo.$

\par\medskip\noindent\textit{Case where $\inf(P_o)=+\infty.$}\ Note
that $\sup(D_o)\geq -\Fo(0)=0>-\infty,$  so that we can apply
Theorem \ref{TLag1}-b. It is enough to prove that
\[
\ls \varphi_0 (0)=+\infty
\]
in the situation where $\varphi_0(0)=\inf(P_o)=+\infty.$ We have
$\ls \varphi_0(0)=\sup_{V\in \mathcal{N}(0)}\inf\{\Fos(\ell);
\ell: \To\ell\in \Co+V\}$ where $\mathcal{N}(0)$ is the set of all
the $\sXY$-open neighbourhoods of $0\in \XXo.$ It follows that for
all $V\in \mathcal{N}(0),$ there exists $\ell\in \LLo$ such that
$\To\ell\in \Co+V$ and $\Fos(\ell)\leq\ls \varphi_0(0).$ This
implies that
\begin{equation}\label{eq-36b}
    \To(\{\Fos\leq\ls\varphi_0(0)\})\cap (\Co+V)\not =\emptyset,\quad
    \forall V\in \mathcal{N}(0).
\end{equation}
On the other hand, $\inf(P_o)=+\infty$ is equivalent to: $
\To(\dom \Fos)\cap \Co=\emptyset. $
\\
Now, we prove ad absurdum that $\ls\varphi_0(0)=+\infty.$ Suppose
that $\ls\varphi_0(0)<+\infty.$ Because of $\To(\dom \Fos)\cap
\Co=\emptyset,$ we have a fortiori
\[
\To(\{\Fos\leq\ls\varphi_0(0)\})\cap \Co=\emptyset.
\]
As $\Fos$ is inf-compact (Lemma \ref{L6}-a) and $\To$ is
continuous (Lemma \ref{L2}-c), $\To(\{\Fos\leq\ls\varphi_0(0)\})$
is a $\sXY$-compact subset of $\XXo.$ Clearly, it is also convex.
But $\Co$ is assumed to be closed and convex, so that by
Hahn-Banach theorem, $\Co$ and $\To(\{\Fos\leq\ls\varphi_0(0)\})$
are \emph{strictly} separated. This contradicts (\ref{eq-36b}),
considering open neighbourhoods $V$ of the origin in
(\ref{eq-36b}) which are open half-spaces. Consequently,
$\ls\varphi_0(0)=+\infty.$ This completes the proof of the
proposition.
\endproof

\subsection{Primal attainment and dual equality}
We are going to consider the following duality diagram, see
Section \ref{sec:problems}:

\begin{equation}
 \begin{array}{ccc}
\Big\langle\ \UU & , & \LL \ \Big\rangle \\
T^\ast  \Big\uparrow & & \Big\downarrow
 T
\\
\Big\langle\ \YY & , & \XX\ \Big\rangle
\end{array}
 \tag{Diagram 1}
\end{equation}
Note that the inclusions $T\LL\subset\XX$ and  $T^\ast
\YY\subset\UU$ which are stated in Lemma \ref{L2} are necessary to
validate this diagram.
\\
Let $F_1, G_1$ and $\gamma_1$ be the  analogues of $F_0,$ $G_0$
and $\gamma_0.$ Denoting $\varphi_1$ the primal value function, we
obtain
\begin{eqnarray*}
F_1(\ell,x) &=& \Fs(\ell)+\iota_{\Ci}(T\ell+x),\quad \ell\in\LL, x\in\XX\\
G_1(y,u)    &=& \inf_{x\in\Ci}\yx -\F(T^\ast y+u),\quad y\in\YY,
u\in\UU\\
\varphi_1(x) &=& \inf\{\Fs(\ell); \ell\in\LL: T\ell\in \Ci-x\},
\quad x\in\XX
\\
\gamma_1(u)&=& \sup_{y\in\YY}\{\inf_{x\in\Ci}\yx -\F(T^\ast
y+u)\},\quad u\in\UU
\end{eqnarray*}
It appears that the primal and dual problems are \Pi\ and \Di.

\begin{lemma}
\label{C4} Assuming \HF\ and \HT, the problems \Po\ and \Pi\ are
equivalent: they have the same solutions and
$\inf(P_o)=\inf(P)\in[0,\infty].$
\end{lemma}
\proof It is a direct consequence of $\dom\Fos\subset\LL,$
$\To\LL\subset\XX$ and $\Fos=\Fs$ on $\LL,$ see Lemma \ref{L2}-a,b
and Lemma \ref{res-05}-c.
\endproof

\begin{proposition}[Primal attainment and dual equality]\label{P3}
Assume that \HF\ and \HT\ hold.
\begin{enumerate}
    \item[(a)] For all $x$ in $\XX,$ we have the little dual equality
        \begin{equation}\label{ped}
            \inf\{\Fos(\ell); \ell\in\LLo, \To\ell=x\}=\Ls(x)\in
            [0,\infty].
        \end{equation}
\end{enumerate}
Assume that in addition \HC\ holds.
\begin{enumerate}
    \item[(b)]  We have the dual equalities
\begin{align}
   & \inf(P_o)=\sup(D)\in[0,\infty]\label{ed1}\\
   & \inf(P_o)=\inf(P)=\inf_{x\in\Ci}\Ls(x)\in [0,\infty]\label{ed1bis}
\end{align}
    \item[(c)] If in addition $\inf(P_o)<\infty,$ then \Po\ is attained in $\LL.$
    \item[(d)] Let $\lh\in\LL$ be a solution to $(P),$ then
$\xh:= T\lh$ is a solution to \PX\ and
$\inf(P)=\Fs(\lh)=\Ls(\xh).$
\end{enumerate}
\end{proposition}

\proof $\bullet$\quad
 We begin with the proof of (\ref{ed1}). As
$\inf(P_o)=\inf(P)$ by Lemma \ref{C4}, we have to show that
$\inf(P)=\sup(D).$ We consider separately the cases where
$\inf(P)<+\infty$ and $\inf(P)=+\infty.$

\par\medskip\noindent\textit{Case where $\inf(P)<+\infty.$}\
Because of \HC, $F_1$ is jointly convex and $F_1(\ell,\cdot)$ is
$\sXYi$-closed convex for all $\ell\in\LL.$ As $T^\ast
\YY\subset\UU$ (Lemma \ref{L2}), one can apply the approach of
Section \ref{sec:convexmin} to the duality Diagram 1. Therefore,
by Theorem \ref{TLag1}-b', the dual equality holds  if $\gamma_1$
is $\sULi$-\usc\ at 0. As in the proof of Proposition \ref{P1}, we
have $-\gamma_1(u)\leq\F(u),$ for all $u\in\UU.$ But $\F$ is the
$\sULi$-\lsc\  regularization of $\Fo+\iota_{\UUo}$ on $\UU$ and
$\Fo$ is bounded above by $1$ on the ball $\{u\in\UUo ;
|u|_\F<1\}.$ As $\LL=(\UU,\NF)',$ $\F$ is also the
$\NF$-regularization of $\Fo+\iota_{\UUo}.$ Therefore, $\F$ is
bounded above by $1$ on $\{u\in\UU ; |u|_\F<1\},$ since
$\{u\in\UUo ; |u|_\F<1\}$ is $\NF$-dense in $\{u\in\UU ;
|u|_\F<1\}.$ As $-\gamma_1 (\leq\F)$ is convex and bounded above
on a $\NF$-neighbourhood of $0,$ it is $\NF$-continuous on
$\icordom (-\gamma_1)\ni 0.$ Hence, it is $\sULi$-\lsc\  at $0.$

\par\medskip\noindent\textit{Case where $\inf(P)=+\infty.$}\
This proof is a transcription of  the second part of the proof of
Proposition \ref{P1}, replacing $\To$ by $T,$ $\Co$ by $\Ci,$ all
the subscripts 0 by 1 and using the preliminary results:
 $\Fs$ is inf-compact (Lemma \ref{L6}) and $T$ is weakly
continuous (Lemma \ref{L2}-h). This completes the proof of
(\ref{ed1}).

\noindent $\bullet$\quad The identity (\ref{ped}) is simply
(\ref{ed1}) with $\Ci=\{x\}.$

\noindent $\bullet$\quad Let us prove (c). By Lemma \ref{L2}-h,
$T$ is $\sLUi$-$\sXYi$-continuous. Since $\Ci$ is $\sXYi$-closed,
$\{\ell\in\LL ; T\ell\in\Ci\}$ is $\sLUi$-closed. As $\Fs$ is
$\sLUi$-inf-compact (Lemma \ref{L6}), it achieves its infimum on
the closed set $\{\ell\in\LL ;T\ell\in\Ci\}$ if
$\inf(P)=\inf(P_o)<\infty.$

\noindent $\bullet$\quad Let us prove (\ref{ed1bis}). The dual
equality (\ref{ed1}) gives us for all $x_o\in\Ci,$
$\inf(P)=\sup_{y\in\YY}\{\inf_{x\in\Ci}\yx-\La(y)\} \leq
\sup_{y\in\YY}\{\langle x_o,y\rangle-\La(y)\}=\Ls(x_o).$ Therefore
\begin{equation}
  \label{eq-03}
  \inf(P)\leq\inf_{x\in\Ci}\Ls(x).
\end{equation}
In particular, equality holds instead of inequality if
$\inf(P)=+\infty.$ Suppose now that $\inf(P)<\infty.$ From
statement (c), we already know that there exists $\lh\in\LL$ such
that $\xh:= T\lh\in\Ci$ and $\inf(P)=\Fs(\lh).$ Clearly
$\inf(P)\leq\inf\{\Fs(\ell) ; T\ell=\xh,
\ell\in\LL\}\leq\Fs(\lh).$ Hence, $\inf(P)=\inf\{\Fs(\ell) ;
T\ell=\xh, \ell\in\LL\}.$ By the little dual equality (\ref{ped})
we have $\inf\{\Fs(\ell) ; T\ell=\xh, \ell\in\LL\}=\Ls(\xh).$
Finally, we have obtained $\inf(P)=\Ls(\xh)$ with $\xh\in\Ci.$
Together with (\ref{eq-03}), this leads us to the desired
identity: $\inf(P)=\inf_{x\in\Ci}\Ls(x).$

\noindent $\bullet$\quad Finally, (d) is a by-product of the proof
of (\ref{ed1bis}).
\endproof

\begin{corollary}\label{res-09}
We have $\dom\Ls\subset\dom\Los,$ $\dom\Ls\subset\XX$ and
    in restriction to $\XX,$ $\Los=\Ls.$
\end{corollary}
\proof The first part is already proved at Lemma \ref{L2}-a.  The
matching $\Los=\Ls$ follows from (\ref{ped0}) and (\ref{ped}).
\endproof

\begin{lemma}\label{L6b}
Under the hypotheses \HF\ and \HT,
 $\Ls$ is $\sXYi$-inf-compact.
\end{lemma}
\proof  By (\ref{ped}): $\inf\{\Fs(\ell); \ell\in\LL,
T\ell=x\}=\Ls(x)$ for all $x\in\XX$ (note that $\Fos=\Fs$ on $\LL$
by Lemma \ref{res-05}-c.) As $T$ is continuous (Lemma \ref{L2}-h)
and $\Fs$ is inf-compact(Lemma \ref{L6}), it follows that $\Ls$ is
also inf-compact.
\endproof

\subsection{Dual attainment}
\label{sec:dualatt} We now consider the following duality diagram
\begin{equation}
 \begin{array}{ccc}
\Big\langle\ \LL & , & \LL^\ast \ \Big\rangle \\
T \Big\downarrow & & \Big\uparrow
 T^\ast
\\
\Big\langle\ \XX & , & \XX^\ast\ \Big\rangle
\end{array}
\tag{Diagram 2}
\end{equation}
where the topologies are the respective weak topologies. The
associated perturbation functions are
\begin{eqnarray*}
  F_2(\ell,x) &=&
  \Fs(\ell)+\iota_{\Ci}(T\ell+x),\quad
  \ell\in\LL, x\in\XX\\
G_2(\zeta,\omega) &=& \inf_{x\in \Ci}\xo-\Fb(T^\ast
\omega+\zeta),\quad \zeta\in\LL^\ast, \omega\in\XX^\ast
\end{eqnarray*}
 As $F_2=F_1,$ the primal problem is \Pi\ and its value function
 is $\varphi_1:$
 \begin{equation}\label{eq-104}
    \varphi_1(x) = \inf_{x'\in \Ci-x}\Ls(x'),\quad x\in\XX
\end{equation}
where we used (\ref{ped}). The dual problem is \Dii.

\begin{proposition}[Dual attainment]\label{res-06a}
Assume that \HF, \HT\ and \HC\ hold. Suppose that
 \begin{equation}\label{eq-106}
    \Co\cap\icordom\Ls\not=\emptyset.
\end{equation}
    Then the dual problem \Dii\ is attained in $\XX^\ast.$
\end{proposition}

\proof
 As $F_2=F_1,$  one can apply the approach of Section
\ref{sec:convexmin} to the duality Diagram 2. Let us denote
$\fss_1$ the $\sXYi$-\lsc\  regularization of $\varphi_1$ and
$\fss_2$ its $\sXYii$-\lsc\  regularization. Since $\XX$ separates
$\YY,$ the inclusion $\YY\subset\XX^\ast$ holds. It follows that
$\fss_1(0)\leq\fss_2(0)\leq\varphi_1(0).$ But we have (\ref{ed1})
which is $\fss_1(0)=\varphi_1(0).$ Therefore, one also obtains
$\fss_2(0)=\varphi_1(0)$ which is the dual equality
\begin{equation}
  \label{ed2}
 \inf(P)=\sup(\overline{D})
\end{equation}
and one  can apply Theorem \ref{TLag1}-c which gives
\begin{equation}\label{eq-114}
     \mathrm{argmax}(\overline{D})=-\partial\varphi_1(0).
\end{equation}
It remains to show that the value function $\varphi_1$  given at
(\ref{eq-104}) is such that
\begin{equation}\label{eq-105}
    \partial\varphi_1(0)\not=\emptyset.
\end{equation}
As the considered dual pairing  $\langle\XX,\XX^\ast\rangle$ is
the saturated algebraic pairing, for (\ref{eq-105}) to be
satisfied, by the geometric version of Hahn-Banach theorem, it is
enough that $0\in\icordom\varphi_1.$ But this holds provided that
the constraint qualification (\ref{eq-106}) is satisfied.
\endproof

Supposing that $\inf(P_o)<\infty$ one knows by Proposition
\ref{P3}-d that \PX\ admits at least a solution $\xh=T\lh$ where
$\lh$ is a solution to \Pi. Let us consider the following new
minimization problem
\begin{equation}
    \textsl{minimize } \Fs(\ell) \quad \textsl{subject to }\quad
T\ell=\xh,\quad \ell\in\LL   \tag{$P^{\xh}$}
\end{equation}
Of course $\lh$ is a solution to \Pi\ if and only if it is a
solution to \Pxh\ where $\xh=T\lh.$ Since our aim is to derive a
representation formula for $\lh,$ it is enough to build our
duality schema upon \Pxh\ rather than upon \Pi. The associated
perturbation functions are
\begin{eqnarray*}
  F_2^{\xh}(\ell,x) &=&
  \Fs(\ell)+\iota_{\{\xh\}}(T\ell+x),\quad
  \ell\in\LL, x\in\XX\\
G_2^{\xh}(\zeta,\omega) &=& \langle\xh,\omega\rangle-\Fb(T^\ast
\omega+\zeta),\quad \zeta\in\LL^\ast, \omega\in\XX^\ast
\end{eqnarray*}
 As $F_2^{\xh}$ is $F_1$ with $\Ci=\{\xh\},$ the primal problem is \Pxh\ and its value
 function is
 \begin{displaymath}
   \varphi_1^{\xh}(x) = \Ls(\xh-x),\quad x\in\XX.
 \end{displaymath}
The dual problem is
\begin{equation}
    \textsl{maximize }\quad \langle\xh,\omega\rangle-\Lb(\omega),
  \quad\omega\in\XX^\ast
  \tag{$\overline{D}^{\xh}$}
\end{equation}

\begin{corollary}[Dual attainment]\label{res-06b}
Assume that \HF\ and \HT\ hold. Suppose that $
\Co\cap\dom\Ls\not=\emptyset.$ Then, $\inf(P_o)<\infty$ and we
know (see
    Proposition \ref{P3}-d) that \PX\
    admits at least a solution. If in addition, there exists a solution
    $\xh$ to \PX\ such that
\begin{equation}\label{CQbis}
    \xh\in\diffdom\Ls,
\end{equation}
then the dual problem \Dxh\ is attained in $\XX^\ast.$
\end{corollary}

\proof Let us specialize Proposition \ref{res-06a} to the special
case where $\Ci=\{\xh\}.$ The dual equality (\ref{ed2}) becomes
\begin{equation}\label{ed3}
 \inf(P^{\xh})=\sup(\overline{D}^{\xh})
\end{equation}
and (\ref{eq-105}) becomes
$\partial\varphi_1^{\xh}(0)\not=\emptyset$  which is implied by
(\ref{CQbis}).
\endproof

\begin{remark}
Let us denote the extended real functions on $\XX^\ast$
    \begin{eqnarray*}
      \La_1&:=&\La+\iota_\YY \\
      \La_2&:=&\Lb
    \end{eqnarray*}
We also denote $\Ls_1,$  $\Ls_2$ their convex conjugates with
respect to $\langle\XX,\XX^\ast\rangle$ and $\widetilde{\La}_1,$
$\widetilde{\La}_2$ their convex $\sYXii$-\lsc\ regularizations.
Clearly, $$\Ls_1=\Ls$$ and the dual equality (\ref{ed3}) is
\begin{equation}\label{eq-112}
    \Ls_1=\Ls_2
\end{equation}
which implies the identity
\begin{equation}\label{eq-111}
\widetilde{\La}_1=\widetilde{\La}_2
\end{equation}
Usual results about convex conjugation tell us that
$\Ls_1(\xh)=\sup_{\omega\in\XX^\ast}\{\langle\xh,\omega\rangle-\widetilde{\La}_1(\omega)\}=\sup(\overline{D}^{\xh})$
and the above supremum is attained at $\ob$ if and only if
$\ob\in\partial_{\XX^\ast}\Ls(\xh).$ This is the attainment
statement in Corollary \ref{res-06b}.
\end{remark}

\subsection{Dual representation of the minimizers}

We keep the framework of Diagram 2 and derive the KKT relations in
this situation. The Lagrangian associated with $F_2=F_1$ and
Diagram 2 is for any $\ell\in\LL, \omega\in\XX^\ast,$
\begin{eqnarray*}
  K_2(\ell, \omega)&:=&\inf_{x\in\XX}\{\xo+\Fs(\ell)+\iota_{\Ci}(T\ell+x)\},\\
&=& \Fs(\ell)-\langle T\ell,\omega\rangle + \inf_{x\in\Ci}\xo.
\end{eqnarray*}

\begin{proposition}[Dual representation]\label{res-07}
Assume that \HF, \HT\ and \HC\ hold.
\\
For any  $\lh\in\LL$ and $\ob\in\XX^\ast,$
 \begin{equation}\label{eq-96}
    \left\{\begin{array}{cl}
      (a) & T\lh\in \Co \\
      (b) & \langle \ob,T\lh\rangle\leq \langle \ob,x\rangle
      \textrm{ for all }x\in \Ci \\
      (c) & \lh\in\partial_{\LL}\Fb(T^*\ob) \\
    \end{array}\right.
\end{equation}
is equivalent to
\begin{equation}\label{eq-98}
    \left\{%
\begin{array}{l}
    \hbox{$\lh$ is a solution to \Po,} \\
    \hbox{$\ob$ is a solution to \Dii\ and} \\
    \hbox{the dual equality (\ref{ed1}) holds.} \\
\end{array}%
\right.
\end{equation}
\end{proposition}

It is well-known that the representation formula (\ref{eq-96}-c):
        \begin{equation}\label{eq-110}
    \lh\in\partial_{\LL}\Fb(T^*\ob)
\end{equation}
 is equivalent to
 \begin{equation*}
    T^*\ob\in\partial_{\LL^*}\Fs(\lh)
\end{equation*}
 and also equivalent to Young's identity
\begin{equation}\label{eq-107}
    \Fs(\lh)+\Fb(T^*\ob)=\langle \ob,T\lh\rangle.
\end{equation}

\proof This proof is an application of Theorem \ref{TLag2}. Under
the general assumptions  \HF, \HT\ and \HC, we have seen at
Proposition \ref{res-06a} that the dual equalities (\ref{ed2}) and
(\ref{ed3}) hold true. Hence, (\ref{eq-98}) is equivalent to
$(\lh,\ob)$ is a saddle-point. All we have to do now is to show
that (\ref{eq-96}) is a translation of the KKT relations
(\ref{KTa}) and (\ref{KTb}).
\\
With $K_2$ as above,  (\ref{KTa}) and (\ref{KTb}) are
$\partial_\ell K_2(\lh,\ob)\ni 0$ and
$\partial_\omega(-K_2)(\lh,\ob)\ni 0.$ Since $-\langle
T\ell,\omega\rangle$ is locally weakly upper bounded as a function
of $\omega$ around $\ob$ and as a function of $\ell$ around $\lh,$
one can apply (Rockafellar, \cite{Roc74}, Theorem 20) to derive
$\partial_\ell K_2(\lh,\ob)=\partial\Fs(\lh)-T^\ast \ob$ and
$\partial_\omega(-K_2)(\lh,\ob)= \partial(-\inf_{x\in\Ci}\langle
x,\cdot\rangle)+T\lh.$ Therefore the KKT relations are
\begin{eqnarray}
  T^\ast \ob&\in &\partial\Fs(\lh)\label{eq-108}\\
 -T\lh &\in & \partial(\iota_{-\Ci}^*)(\ob)\label{eq-109}
\end{eqnarray}
where $\iota_{-\Ci}^*$ is the convex conjugate of the convex
indicator of $-\Ci.$
\\
As a convex conjugate, $\Fs$ is a closed convex functions. Its
convex conjugate is $\Fb.$ Therefore (\ref{eq-108}) is equivalent
to the following equivalent statements
\begin{eqnarray*}
  &&\lh \in \partial\Fb(T^\ast \ob)\\
&& \Fs(\lh)+ \Fb(T^\ast \ob)=\langle \lh,T^*\ob\rangle
\end{eqnarray*}
Similarly, as a convex conjugate $\iota_{-\Ci}^*$ is a closed
convex functions. Its convex conjugate is $\iota_{-\bar\Ci}$ where
$\bar\Ci$ stands for the $\sXYii$-closure of $\Ci.$ Of course, as
$\Ci$ is $\sXYi$-closed by hypothesis \HC, it is a fortiori
$\sXYii$-closed, so that $\bar\Ci=\Ci.$ Therefore (\ref{eq-109})
is equivalent to
\begin{equation}\label{eq-13}
   \iota_{\Ci}(T\lh)+\iota_{-\Ci}^*(\ob)=\langle
-T\lh,\ob\rangle.
\end{equation}
It follows from (\ref{eq-13}) that $\iota_{\Ci}(T\lh)<\infty$
which is equivalent to $ T\lh\in\Ci.$
\\
Now (\ref{eq-13}) is $-\langle
T\lh,\ob\rangle=\iota_{-\Ci}^*(\ob)=-\inf_{x\in\Ci}\langle
x,\ob\rangle$ which is $\langle T\lh,\ob\rangle
=\inf_{x\in\Ci}\langle x,\ob\rangle.$ This completes the proof.
\endproof

\begin{remark}
Thanks to Proposition \ref{P3}-d, (\ref{eq-107}) leads us to
\begin{equation}\label{eq-113}
    \Ls(\xh)+\Lb(\ob)=\langle\xh,\ob\rangle
\end{equation}
for all $\xh\in\dom\Ls$ and all $\ob\in\XX^\ast$ solution to \Dxh.
By Young's inequality:
$\Ls_2(\xh)+\widetilde{\La}_2(\ob)\geq\langle\xh,\ob\rangle$ and
the identities (\ref{eq-112}) and (\ref{eq-113}), we see that
$\widetilde{\La}_2(\ob)\geq\Lb(\ob).$ But, the reversed inequality
always holds true. Therefore, we have $
    \widetilde{\La}_2(\ob)=\Lb(\ob).
$ This proves that
    $
    \Lb=\widetilde{\La}_2 \textrm{ on } \dom\Lb:
    $
$\Lb$ is $\sigma(\XX^*,\XX)$-\lsc\ on its effective domain.
\end{remark}

\begin{proposition}\label{res-08}
Assume that \HF, \HT\ and \HC\ hold. Any solution $\ob$ of \Dii\
or \Dxh\ shares the following properties
\begin{itemize}
 \item[(a)] $\ob$ is in the $\sigma(\XX^*,\XX)$-closure of $\dom\La;$
 \item[(b)] $T^\ast \ob$ is in the $\sigma(\LL^*,\LL)$-closures of
 $T^\ast(\dom\La)$ and $\dom\Fo.$
\end{itemize}
If in addition the level sets of $\F$ are $\NF$-bounded, then
\begin{itemize}
 \item[(a')] $\ob$ is in $\YY''.$ More precisely, it is in the $\sigma(\YY'',\XX)$-closure of $\dom\La;$
 \item[(b')] $T^\ast \ob$ is in $\UU''.$ More precisely, it is in the $\sigma(\UU'',\LL)$-closures of
 $T^\ast(\dom\La)$ and $\dom\Fo$
\end{itemize}
where $\YY''$ and $\UU''$ are the \emph{topological} bidual spaces
of $\YY$ and $\UU.$ This occurs if $\F,$ and therefore $\Fs,$ is
an even function.
\end{proposition}
\proof
 \boulette{a}
Because of (\ref{eq-113}), we have $\ob\in\dom\Lb.$ As
$\widetilde{\La}_2\leq \Lb$ and
$\widetilde{\La}_1=\widetilde{\La}_2$ (see (\ref{eq-111})), we
obtain $\ob\in\dom\widetilde{\La}_1$ which implies that $\ob$ is
in the $\sYXii$-closure of $\dom\La.$

 \Boulette{b}
By Lemma \ref{L2}-d,  $T^*$  is continuous from $\XX^*$ to
$\LL^*.$ It follows from (a) that $T^\ast \ob$ is in the
$\sigma(\LL^*,\LL)$-closure of $T^\ast(\dom\La).$
\\
 On the other hand, $T^\ast \ob\in\dom\Fb$ and
 $\Fb$ is the $\sULii$-closed convex regularization of $\Fo.$ It follows
 that $T^\ast \ob$ is in the $\sigma(\LL^*,\LL)$-closure of $\dom\Fo.$

 \Boulette{a'}
Because of (a), $\ob$ is the $\sigma(\XX^*,\XX)$-limit of a
generalized sequence $\{y_\alpha\}$ in $\dom\La.$ Our additional
assumption allows us to take $\{y_\alpha\}$ in a $\NF$-ball: it is
an equicontinuous set. It follows with \cite[Cor. of Prop.
III.5]{Bou-EVT} that $\ob$ is continuous on $\XX.$

 \Boulette{b'} Similar to (b)'s proof using (a') and Lemma \ref{L2}-f.
\endproof


\end{document}